\documentclass[ejs]{imsart}

\RequirePackage{amsthm,amsmath}
\RequirePackage[numbers]{natbib}
\RequirePackage[colorlinks,citecolor=blue,urlcolor=blue]{hyperref}

\usepackage{amsmath,amsthm,amsfonts,amssymb}
\usepackage{xcolor}
\usepackage{subcaption} 
\usepackage{graphicx}
\usepackage{pgfplots}
\usepackage{enumitem}
\usepackage[ruled,vlined]{algorithm2e}

\arxiv{}

\startlocaldefs

\newcommand{\RR}{\mathbb{R}}
\newcommand{\qv}[1]{\langle  #1 \rangle}
\newcommand{\ip}[1]{\langle  #1 \rangle}
\newcommand{\given}{\,|\,}
\newcommand{\HHH}{\mathbb{H}}
\newcommand{\eps}{\varepsilon}
\newcommand{\FF}{\mathcal{F}}

\newcommand{\gaatp}{\overset{P}{\to}}
\newcommand{\lle}{\lesssim}

\newcommand{\NN}{\mathbb{N}}
\newcommand{\GG}{\mathcal{G}}
\newcommand{\EE}{\mathbb{E}}
\renewcommand{\phi}{\varphi}
\newcommand{\sign}{\mathop{\rm sign}\nolimits}
\def\q{\theta}
\def\g{\gamma}
\def\e{\epsilon}
\def\d{\delta}
\def\a{\alpha}
\def\b{\beta}
\def\s{\sigma}
\def\ra{\rightarrow}
\def\var{\mathop{\rm var}\nolimits}

\theoremstyle{plain}
\newtheorem{thm}{Theorem}[section]
\newtheorem{lemma}[thm]{Lemma}
\newtheorem{cor}[thm]{Corollary}
\newtheorem{prop}[thm]{Proposition}

\newtheorem{asspt}[thm]{Assumption}

\theoremstyle{definition}
\newtheorem{ex}[thm]{Example}
\newenvironment{prf}{\vspace{0ex}\begin{proof}[\bf Proof]}{\end{proof}\vspace{2ex}}

\endlocaldefs

\begin{document}

\begin{frontmatter}

\title{Semi-parametric Bernstein-von Mises Theorem in Linear Inverse Problems}
\runtitle{Semi-parametric BvM in Linear Inverse Problems}


\begin{aug}
	\author[A]{\fnms{Adel} \snm{Magra}\ead[label=e1]{a.magra@vu.nl}}

	\author[B]{\fnms{Aad} \snm{van der Vaart}\ead[label=e2]{a.w.vandervaart@tudelft.nl}}

	\author[A]{\fnms{Harry} \snm{van Zanten}\ead[label=e3]{j.h.van.zanten@vu.nl}}

	\runauthor{Magra, van der Vaart, van Zanten}
	\address[A]{VU Amsterdam}\printead[presep={\ }]{e1,e3}
	
	\address[B]{TU Delft}\printead[presep={\ }]{e2}
	
\end{aug}

\begin{abstract}
We consider a Bayesian approach for the recovery of scalar parameters arising in inverse problems. We consider a general signal-in white noise model where we have access to two independent noisy observations of a function, and of a linear transformation of the function. The linear operator is unknown up to a scalar parameter. We present a Bernstein-von Mises theorem for the marginal posterior of the scalar under regularity assumptions of the operator. We further derive Bernstein-von Mises results for different priors and apply them to  two concrete examples: the recovery of the thermal diffusivity in a heat equation problem, and the recovery of a location parameter in a semi-blind deconvolution problem.

\end{abstract}



\end{frontmatter}


\section{Introduction}

In nonparametric statistical inverse problems the goal is to 
recover a function $f$ from a noisy observation of a transformed 
version $K\!f$ of the function, where $K$ is some known operator, typically
with unbounded inverse. Many methods have been proposed and 
studied for this problem, including  frequentist regularization methods  and 
nonparametric Bayes approaches. 
The focus in the theoretical literature is mostly on the setting 
that the operator $K$ is known and the function $f$ is the unknown object of interest
that needs to be recovered from the data. In several important applications however, 
the operator $K = K_\theta$  depends on an unknown, Euclidean parameter $\theta$, 
and it is that parameter which is actually  the main object of interest. 

An example from biology 
occurs in the paper \cite{gao}, which deals with the modeling of biochemical 
interaction networks. In that case $K= K_\theta$ is the 
  solution operator of a linear ordinary differential equation describing the time evolution of gene expression levels. The function $f$ describes how the activity of a so-called transcription factor changes over time and 
the parameter vector $\theta$  describes important aspects of the chemical reaction that is modeled,
 like the basal transcription rate and the rate of decay of mRNA. 
Another setting in which problems of this type arise naturally is in cosmology. In general relativity for instance, 
Einstein's equations describe how the large-scale structure of the universe evolves 
from small  energy density fluctuations in the early universe. These equations are 
partial differential equations  that depend on a number of important cosmological constants that cosmologists 
want to learn. Given the early state $f$ of the universe, general relativity gives ``forward'' predictions 
$K_\theta f$ of certain observable quantities. By comparing these predictions with (noisy) observations, 
inference can be made about the parameters $\theta$.  
 In such applications the operator $K_\theta$ 
is typically an   operator giving an appropriate solution to Einstein's equations. 
A particular example of this approach occurs in weak lensing, which is the state-of-the-art
method used in  recent and future cosmological surveys, see for instance \cite{lensing}.
In that case the parameter vector $\theta$ of interest includes important constants like 
the matter density of the universe and the matter inhomogeneity.

The natural Bayesian approach to this type of statistical inverse problems with parameter-dependent 
operators, often followed in practice, starts with 
endowing the pair $(\theta, f)$ 
with a prior distribution. The particular data generating mechanism gives rise  
to a likelihood and together the prior and the likelihood yield a posterior distribution for the 
unknown pair $(\theta, f)$. The marginal posterior for $\theta$ can then be used to 
make inference about $\theta$. 
This paper is motivated by the fact that this type of Bayesian methodology is commonly used in inverse problems with parameter-dependent
operators, yet there are no readily applicable theoretical results giving insight into their fundamental performance. 

To obtain some first insight   we study these matters in this paper 
in the context of semi-parametric inverse problems that are interesting yet relatively tractable
mathematically. Let $\Theta$ be a compact subset of $\mathbb{R}$ and let $(H,\ip{\cdot,\cdot},\|\cdot\|)$ be a separable Hilbert space. Assume $f\in H$ and $\forall \theta\in\Theta, \ K_\theta:H \to H$. We consider an observation model consisting of two independent noisy observations, one observation for $f$, and one for $K_\theta f$. Namely
\begin{align}
	\label{eq: obs1} X^{(1)}&= f + \frac{1}{\sqrt{n}}\dot{W}^{(1)},\\
	\label{eq: obs2} X^{(2)}&= K_\theta f + \frac{1}{\sqrt{n}}\dot{W}^{(2)},
\end{align}
where $\dot{W}^{(1)}$ and $\dot W^{(2)}$ are independent iso-Gaussian processes for $H$ and $n\in\mathbb{N}$ is the signal-to-noise 
ratio. An iso-Gaussian process $\dot W$ for $H$ is a stochastic process $\dot W=\{\dot W_h:h\in H\}$ with $\dot W_h\sim N(0,\|h\|^2)$.
(In Section~\ref{sec:result} we describe the model in the statistically equivalent sequence form, 
while in Section~\ref{sec:applications} we use diffusion equations.)

We take a Bayesian approach and endow the pair $(\theta, f)$ 
with a prior distribution by putting a prior on $f$ and 
an independent prior with a positive, continuous Lebesgue density on $\theta$. Assuming that there exists a true pair of parameters $(\theta_0,f_0)$, we present results that describe the behavior of the marginal posterior
for the parameter $\theta$ of interest as $n\to \infty$. 

Our main result is Theorem \ref{thm: main}, a semi-parametric Bernstein-von Mises
(BvM) theorem for the marginal posterior of $\theta$.  It 
gives conditions under which  the marginal posterior behaves asymptotically like a normal distribution 
centered at an efficient estimator of $\theta$, with a variance equal to the inverse of the 
efficient Fisher information. Such a result guarantees in particular the validity of uncertainty 
quantification, in the sense that that credible 
intervals for the parameter $\theta$ obtained from the marginal posterior are also 
asymptotic frequentist confidence intervals in  that case, see for instance  \cite{castillo}.

We assume that the dependence $\theta\mapsto K_\theta$ is smooth, which allows
a semi-parametric score calculus, resulting in a ``least favorable submodel''.
Theorem~\ref{thm: main} then has two ingredients, both involving 
the ``least favorable direction''. The first is a ``bias term'' that arises in an expansion of
the log likelihood in the least favorable direction, and must be suitably small. The second ingredient is 
insensitivity of the prior relative to location shifts in the least favorable direction. These conditions are 
similar to the ones in \cite{castillo}, but take a special, attractive form in our setting.

Verification of the general conditions may involve a posterior contraction rate. This can
be obtained under additional assumptions on the operators $K_\theta$ and the usual conditions 
on the prior on $f$: the existence of sieves with relatively small entropy and enough 
prior mass on neighborhoods of $f_0$, as expressed in Lemma~\ref{thm: main2}.
Remarkably, some examples can be handled by special properties of the operators, without
needing a contraction rate. Corollaries~\ref{cor: rateCondition} and~\ref{cor: constantNorm} 
summarize the interplay between the least favorable direction and the prior.

We apply our general theorem to operators in two different problems: the heat equation and semi-blind deconvolution,
and consider these with two different types of priors: Gaussian and $p$-exponential priors.

Our paper addresses a parameter that is included in an inverse problem in a nonlinear manner. There are many 
other papers that deal with Bernstein-von Mises theorems in inverse problems, often of linear functionals
of a nonparametric parameter, e.g.\ \cite{Knapiketal2011}, \cite{Ray2017}, \cite{ismaeljudith}, \cite{RayvdV}, \cite{nickl2018}, 
\cite{GiordanoKekkonen}, \cite{MonardNicklPternain2021}, \cite{Nickl23}.

\subsection{Organization of the paper}
We present the semi-parametric setting and a general BvM theorem for our setting, along with two corollaries
 in Section~\ref{sec:result}. Next we give example applications in Section~\ref{sec:applications}, and
specific priors in Section~\ref{sec:priors}. In Section~\ref{sec:simulations}, we numerically  illustrate our results, 
using a Metropolis-Hastings scheme for sampling 
from the marginal posterior for $\theta$. Some proofs are given in Section~\ref{sec: proof}.

\subsection{Notation}
The space of square integrable functions on the interval $[0,1]$ is denoted by $L^2[0,1]$, and the space of Hölder continuous functions of order $\alpha$ by $C^\alpha[0,1]$. For $f,g\in L^2[0,1]$ the $L^2$-norm of $f$ is 
$\|f\|^2:=\int_0^1 f(x)^2\,\,\text{d}x$ and the $L^2$-inner product is $\qv{f, g} = \int_0^1 f(x)g(x)\,\text{d}x$.
The norm of an element $h\in H$ of our basic Hilbert space is denoted by $\|h\|$, or occasionally by $\|h\|_H$ for emphasis on $H$.
 For a bounded operator $A: H\to S$ between two normed spaces
$\|A\|_{H \to S} = \sup_{\|h\| \le 1} \|Ah\|_S$ denotes the operator norm, which is abbreviated to $\|A\|$ if $H=S$.
For two numbers $a$ and $b$, we denote by $a\wedge b$ the minimum of $a$ and $b$, and by  $a\vee b$ their maximum. For two sequences $a_n$ and $b_n$, we mean by $a_n\asymp b_n$ that $|a_n/b_n|$ is bounded away from zero and infinity as $n\to\infty$, and by $a_n \lle b_n$ that $a_n/b_n$ is bounded. The total variation distance between two probability measures $P$ and $Q$ defined on the same probability space $(\Omega,\mathcal{F})$ is given by $\|P-Q\|_{\text{TV}}:=2\sup_{F\in\mathcal{F}} |P(F)-Q(F)|$. Finally, for a metric space $(A, d)$ and $\eps > 0$, we denote by  $N(\eps, A, d)$ the minimum number of balls of radius $\eps$ needed to cover $A$.

\section{General results}
\label{sec:result}
We assume that the pair $(\theta, f)$ belongs to $\Theta \times H$, where $\Theta\subset \mathbb{R}$ is compact. We endow the unknown pair with a prior $\Pi$ of the form 
$\Pi = \pi_\theta \times \pi_f$, 
where $\pi_\theta$ has a  continuous Lebesgue density that is bounded away from $0$ and $\infty$ on $\Theta$ and $\pi_f$ is a 
prior on $H$. Let $\{e_k\}_{k\in\mathbb{N}}$ be an arbitrary orthonormal basis for $H$, and let $P^n_{\theta, f}$ be the law of the pair 
$X^n:=(X^{(1)}, X^{(2)})$ 
of sequences $X^{(i)}=(X^{(i)}_1,X^{(i)}_2,\ldots)$ defined by 
\begin{align*}
	 X^{1)}_k&:=\ip{f, e_k} + n^{-1/2}Z_{1,k},\\
	 X^{(2)}_k&:=\ip{K_\theta f, e_k} + n^{-1/2}Z_{2,k},
\end{align*}
where the $Z_{i,k}$ are i.i.d standard Gaussian random variables, for $i=1,2$ and $k=1,2,\ldots$. This model is statistically equivalent
to the white noise model \eqref{eq: obs1}--\eqref{eq: obs2}, in the sense that the likelihood ratios are equal.  We shall also write $X^{(i)}_k=\ip{X^{(i)}, e_k}$
and  for $h\in H$, more generally $\ip{X^{(i)},h}:=\sum_k X^{(i)}_k\ip{h,e_k}$. The latter 
random series is convergent, both in $L_2$ and almost surely.

Let $P_0^n$ be the distribution of $X^n=(X^{(1)}, X^{(2)})$ if $f=0$, or 
equivalently the distribution of the noise process $(n^{-1/2}\dot W^{(1)}, n^{-1/2}W^{(2)})$. 
The log-likelihood has the following expression (e.g.\ Lemma~L.4 in \cite{vaartghosal}),
\begin{align}
	\log \frac{dP^n_{\theta,f}}{dP^n_0}(X^{(1)},X^{(2)}) &= n\ip{X^{(1)},f} - \frac{n}{2}\|f\|^2
	+ n\ip{X^{(2)},K_\theta f}-\frac{n}{2}\|K_\theta f\|^2.
\label{EqLikelihood} 
	\end{align}
Combined with the prior this results in a posterior distribution $\Pi(\cdot \given X^n)$ for the pair $(\theta, f)$. 
We are in particular interested in the marginal posterior $B \mapsto \Pi(\theta \in B\given X^n)$.

Since  the parameter of interest $\theta$ is real-valued and $f$ is an infinite-dimensional nuisance parameter, 
the behavior of this marginal posterior is determined by the semi-parametric structure of the 
model. Because the model as given is equivalent to observing $n$ i.i.d.\ copies of the same model but with $n=1$, the usual
semi-parametric theory (as in \cite{vdV91}, \cite{BKRW}, or Chapter~25 of \cite{AadsBoek}) applies. Alternatively,
we may use the fact that the model is (exactly) locally asymptotically normal and apply the general theory given in \cite{McNeney}. 
In Section~\ref{sec: proof} ahead we apply the Bayesian semi-parametric framework of Chapter~12 in \cite{vaartghosal},
which is based on \cite{castillo}.

The key is local asymptotic normality (LAN), an expansion of the log  likelihood in the
neighborhood of a fixed parameter $(\theta,f)$. This is equivalent to ordinary log likelihood
expansions along one-dimensional submodels indexed by parameters of the form $t\mapsto (\theta+t,f+tg)$, for
$t\in\RR$ and fixed ``directions''   $g\in H$. We assume that the operators $K_\theta: H\to H$ are differentiable
in $\theta$ in the sense that there exist operators $\dot K_\theta: H\to H$ such that, for $f\in H$, as $s \to 0$,
\begin{equation}
\label{EqDerivativeK}
\frac{K_{\theta+s} f - K_{\theta} f}{s} \to \dot K_\theta f,\qquad \text{ in }H.
\end{equation}
Then we have the following lemma. Let $K_\theta^*: H\to H$ be the adjoint of $K_\theta$.

\begin{lemma}\label{lem:LAN}
If the preceding display holds, then for every $g\in H$ and $n\rightarrow \infty$, 
with $\dot W^{(1)}$ and $\dot W^{(2)}$ defined in \eqref{eq: obs1}--\eqref{eq: obs2}, in $P_{\q,f}^n$-probability,
	\begin{equation}\label{eq: lan}
		\begin{split}
			\log \frac{dP_{\theta+t/\sqrt n, f+tg/\sqrt n}^n}{dP_{\theta, f}^n} &\to
			t\ip{g,\dot W^{(1)}} + t\ip{\dot K_{\theta}f+K_{\theta}g,\dot W^{(2)}}\\
			& \qquad\qquad -\frac{1}2t^2\|g\|^2 -\frac{1}2t^2\|\dot K_\theta f + K_\theta g\|^2.
		\end{split}
	\end{equation}
Furthermore, the Fisher information $\|g\|^2 + \|\dot K_\theta f + K_\theta g\|^2$ is minimized over $g$ at $g=-\gamma_{\theta,f}$ given by
\begin{equation}
  \gamma_{\theta, f} = (I+K^*_\theta K_\theta)^{-1}K^*_\theta\dot K_\theta f. 
\label{EqLeastFavourableDirection}
\end{equation}
If $\dot K_\q f\not=0$, then the minimal value $\tilde I_{\q,f}:=\|\g_{\q,f}\|^2+\|\dot K_\theta f - K_\theta \g_{\q,f}\|^2$ is strictly positive.
\end{lemma}

The perturbation $g=0$ corresponds to the model in which $f$ is known, in which case the ``Fisher information''
in the lemma reduces to $I_\theta:=\|\dot K_\theta f \|^2$, and is the ordinary Fisher information 
for $\theta$ (or for $\theta+t$ at $t=0$) in the parametric model
with $\theta$ as the only parameter. This information is positive if and only if $\dot K_\q f\not=0$.
For a general perturbation $g$ the ``Fisher information'' for $t$ in the lemma is the information 
in the submodel $t\mapsto (\theta+t, f+tg)$.
The submodel with $g=-\gamma_{\theta,f}$ provides the smallest information and in this sense
is ``least favorable'' for estimating $\theta$.

It is insightful to obtain these quantities also in terms of score functions, as is usual in semi-parametric theory 
(see e.g.\ \cite{vdV91}, \cite{McNeney}). The score function for $\theta$ in the model with only $\theta$ as a parameter
(i.e.\ $g=0$) is equal to $\ip{\dot K_{\theta}f,\dot W^{(2)}}$. The other part of the linear term in the expansion in Lemma~\ref{lem:LAN}
is the score function for the ``nuisance parameter'' $f$ in the direction of
$g$, given by $\ip{g,\dot W^{(1)}} + \ip{K_{\theta}g,\dot W^{(2)}}$. By definition the \emph{efficient score function}
for $\theta$ is the score function for $\theta$ minus its projection onto the closed linear span of the scores for
the nuisance parameter $f$. This leads to minimizing the square distance
$$g\mapsto \mathbb{E}_{\q,f}\bigl[\ip{\dot K_{\theta}f,\dot W^{(2)}}-\ip{g,\dot W^{(1)}} - \ip{K_{\theta}g,\dot W^{(2)}}\bigr]^2
=\|g\|^2+\|\dot K_\theta f-K_\theta g\|^2.$$
The solution is the \emph{least favorable direction} $g=\gamma_{\theta,f}$ given in \eqref{EqLeastFavourableDirection}, 
and the minimum value is the \emph{efficient
Fisher information} $\tilde I_{\theta, f}$. The latter can be obtained by substituting $g=\gamma_{\theta,f}$, and also as
$\|\dot K_\theta f\|^2- \|\gamma_{\theta, f}\|^2-\|K_\theta\gamma_{\theta,f}\|^2$, by Pythagoras's rule and the orthogonality
of projection. This can be written as $\|\dot K_\theta f\|^2- \ip{\gamma_{\theta, f},(I+K_\q^*K_\theta)\gamma_{\theta,f}}$,
or in view of \eqref{EqLeastFavourableDirection},
\[
\tilde I_{\theta, f} = \|\dot K_\theta f\|^2 - \ip{K^*_\theta\dot K_\theta f, 
	(I+K^*_\theta K_\theta)^{-1}K^*_\theta\dot K_\theta f}.
\]
The expression shows that $\tilde I_{\q,f}$  is strictly less than the ordinary Fisher information $I_\theta$
if $\gamma_{\q,f}\not=0$, in which case there is a loss of information due to the fact that $f$ is unknown.

In this setting we say that the \emph{Bernstein-von Mises theorem holds at $(\q_0,f_0)$} if,
in $P_{\theta_0, f_0}^n$-probability, as $n\to\infty$, 
  \begin{align}
         \Big\|\Pi(\theta \in \cdot \given X^n) - 
         N\Big(\theta_0 + \frac{1}{\sqrt{n}} \Delta_{ \theta_0, f_0}^n, 
        \frac{1}{n}\tilde I^{-1}_{\theta_0, f_0}\Big)\Big\|_{\text{TV}}\to 0,
\label{EqBvMTheorem}
    \end{align}
where  $\Delta_{ \theta_0, f_0 }^n$ are measurable
transformations of $X^n$ such that $\Delta_{ \theta_0, f_0 }^n\sim N(0, \tilde I^{-1}_{\theta_0, f_0})$. 
In all our results the latter variables satisfy
$$\Delta_{\q,f}^n
=\frac1{\tilde I_{\q,f}}\bigl(\ip{\dot K_{\theta}f,\dot W^{(2)}}-\ip{\gamma_{\theta,f},\dot W^{(1)}} - \ip{K_{\theta}\gamma_{\theta,f},\dot W^{(2)}}\bigr).$$
In other words, the variables $\Delta_{ \theta_0, f_0 }^n$ are the efficient score functions at $(\q_0,f_0)$ divided by
the efficient Fisher information. The Bernstein-von Mises theorem \eqref{EqBvMTheorem} uses the
total variation norm, which is much stronger than the bounded Lipschitz metric used more often
recently (e.g.\ \cite{ismaeljudith},\cite{RayvdV}).

We prove the Bernstein-von Mises theorem under general conditions. 
We assume that the operator $K_\q$ is ``smoothing'' in that its range $K_\q H$ belongs to a
normed space $S\subset H$ with a stronger norm $\|\cdot\|_S$ than $H$.
For instance, if $H$ is a space of Lebesgue square integrable functions, then $S$ may be 
a Sobolev type space of a certain order, as we shall see in the examples of Section~\ref{sec:applications}. 
We assume that the norm of $S$ is strong enough so that its unit ball has a finite entropy integral in $H$:
\begin{align}\label{eq: integralCondition} 
\int_0^1 \sqrt{\log N\bigl(\eps, \{h: \|h\|_S\le 1\}, \|\cdot\|_H\bigr)}\,d\eps <\infty.
\end{align}
Furthermore, we assume that the map $\theta\mapsto K_\theta$ satisfies the following regularity conditions.

\begin{asspt}\label{asspt:main}
The derivative in \eqref{EqDerivativeK} exists and there exist constants $D_0, D_1,D_2,D_3,D_4$ such that  $\forall \theta, \theta_1, \theta_2\in\Theta$ and $\forall f_1, f_2 \in H$, 
\begin{enumerate}[label=(\roman*),leftmargin=*, itemsep=0.5ex, before={\everymath{\displaystyle}}]%
\item $\|K_\theta\| \le D_0, \qquad \|\dot K_\theta\| \le D_0$,
\item $	\|{K_\theta - K_{\theta_0}}\|_{H\to S} \le  	D_1|\theta - \theta_0|$,
\item $	\bigl\|K_\theta - K_{\theta_0} -(\theta-\theta_0)\dot K_{\theta_0}\bigr\|_{H \to S} \le D_2 (\theta - \theta_0)^2$,
\item $ \|K_{\theta_1} f_1 - K_{\theta_2}f_2\|   \le  D_3 |\theta_1-\theta_2| (\|f_1\| \wedge \|f_2\|)+D_3\|f_1-f_2\|$,
\end{enumerate}
\end{asspt}

Our main result asserts that the semi-parametric BvM theorem holds 
for a  ``good enough'' prior $\pi_f$ on $f$. Here ``good enough'' means that the posterior 
of $(\q,f)$ contracts around the true $(\theta_0,f_0)$ and is insensitive to shifts in the least favorable direction. 

We  denote the true parameter pair by $(\theta_0, f_0)$, and write
$\Pi^{\theta=\theta_0}(\cdot\given X^n)$ for  the posterior distribution of $f$ in the model where $\theta_0$ is known.
The posterior of $(\q,f)$ is said to be \emph{consistent} at  $(\q_0,f_0)$ if 
$\Pi\bigl(\|\q-\q_0\|<\e_n,\|f-f_0\|<\e_n| X^n)\to 1$ in $P_{\q_0,f_0}^n$ probability, 
for some sequence $\e_n\ra0$, called the \emph{rate of contraction}.

\begin{thm}\label{thm: main}
Assume that \eqref{eq: integralCondition}  and Assumption~\ref{asspt:main}(i)-(iii) are satisfied, that  the posterior distribution of
$(\q,f)$ is consistent at $(\q_0,f_0)$, and that $\dot K_{\q_0} f_0\not=0$. Furthermore, 
assume that there exists a sequence $\{\gamma_n\}\subset H$ with $\g_n\to\g_{\q_0,f_0}$
and subsets $\Theta_n\subset\Theta$ and $H_n \subset H$ such that
$\sqrt n(\Theta_n-\q_0)\uparrow\RR$ and 
\begin{align}\label{eq: contraction}
\Pi(\Theta_n\times H_n | X^n)\gaatp 1, 
\qquad \inf_{\theta\in\Theta_n}\Pi^{\theta=\theta_0}\bigl(H_n + (\theta-\theta_0)\gamma_n | X^n\bigr) &\gaatp 1,\\
\label{eq: lfdrate}
	\sqrt{n}\sup_{f\in H_n}\bigl|\ip{f-f_0,(I+K^*_{\theta_0}K_{\theta_0})(\gamma_{\theta_0,f}-\gamma_n)}\bigr|&\to 0,\\
\label{eq: priorshift}
	\sup_{\theta\in\Theta_n, f\in H_n} \left | \frac{\log(d\pi_{f+(\theta-\theta_0)\gamma_n}/d\pi_f)(f)}{1+n(\theta-\theta_0)^2}\right | &\to 0.
\end{align}
Then the Bernstein-von Mises theorem \eqref{EqBvMTheorem} holds at $(\q_0,f_0)$.
 \end{thm}

The sets $H_n$ in the theorem must be large enough to possess posterior mass tending to one \eqref{eq: contraction}, and be small enough so that
\eqref{eq: lfdrate}-\eqref{eq: priorshift} hold. The assumed posterior consistency implies that the sets can always be restricted
to shrinking balls around $f_0$. However, the rate $\sqrt n$ in \eqref{eq: lfdrate} may require a rate of posterior contraction,
and \eqref{eq: priorshift} may depend on still other properties of $H_n$, in interaction with the choice of $\gamma_n$.

Condition \eqref{eq: lfdrate} captures a remainder term in a LAN expansion of the log-likelihood ratio in the least favorable submodel
(analogous to (12.13) in \cite{vaartghosal}). This condition requires a rate or structural properties and is left as a main condition, 
while other, lower order remainder terms are negligible under Assumption~\ref{asspt:main}. 
Condition \eqref{eq: priorshift} 
requires insensitivity of the prior on $f$ under shifts in the least favorable direction. 
Below we replace the condition by the sufficient condition that there exist positive numbers $\e_n\ra 0$ and $\eta_n\ra0$ such that
\begin{align}
\pi_f\bigl(\|f-f_0\|\le \eps_n/(2D_3+2)\bigr) &\ge e^{-n\eps_n^2/65},\label{eq: thmPrior}\\
  \pi_f\Bigl(\exists s, t\in (-\e_n,\e_n): \Bigl| \log\frac{d\pi_{f+(s+t)\gamma_n}} {d\pi_{f+t\g_n}}(f)\Bigr|
> \eta_n (1+nt^2)\Bigr)&\le e^{-3n\e_n^2}.\label{EqPriorLikelihoodRatio}
\end{align}
We verify these conditions for several classes of priors in  Section~\ref{sec:priors}.

The choice of $\gamma_n$ is pivotal to verifying conditions \eqref{eq: lfdrate} and \eqref{eq: priorshift}
simultaneously. Selecting $\gamma_n = \gamma_{\theta_0,f}$ renders condition \eqref{eq: lfdrate} trivial,
but the dependence on $f$ (which would be permitted) may make \eqref{eq: priorshift} difficult to check. Selecting 
$\g_n=\g_{\q_0,f_0}$ reduces (in view of \eqref{EqLeastFavourableDirection}) condition \eqref{eq: lfdrate}  to
\begin{align}\label{eq: lfdrate2}
	\sqrt{n}\sup_{f\in H_n} \bigl|\ip{K_{\theta_0}(f-f_0),\dot K_{\theta_0}(f-f_0)}\bigr| \to 0.
\end{align}
This condition may be valid due to special properties of the operators, or may require a posterior rate of contraction for $f$.
For a general choice $\gamma_n \to \gamma_{\theta_0,f_0}$, the Cauchy-Schwarz inequality and Assumption~\ref{asspt:main}
allow to reduce \eqref{eq: lfdrate} to
\begin{equation}
\sqrt{n}\sup_{f\in H_n} \|f-f_0\|\, \|\g_{\q_0,f}-\g_n\|\to 0.
\label{eq: lfdrate3}
\end{equation}
Here $\g_{\q_0,f_0}$ could replace $\g_{\q_0,f}$, with the difference bounded by $\|\g_{\q_0,f}-\g_{\q_0,f_0}\|\lesssim \|f-f_0\|$, 
in view of \eqref{EqLeastFavourableDirection} and Assumption~\ref{asspt:main}.
These observations lead to the following corollary.

\begin{cor}\label{cor: rateCondition}
Let \eqref{eq: integralCondition} and Assumption~\ref{asspt:main}  be satisfied, and
assume that $\dot K_{\q_0} f_0\not=0$ and that $\q_0$ is interior to $\Theta$. 
Assume that the posterior distributions of $(\q,f)$ and of $f$ given $\q_0$ contract at rate $\e_n\gg n^{-1/2}$.
Suppose that \eqref{eq: thmPrior} and \eqref{EqPriorLikelihoodRatio} are  satisfied for some
$\gamma_n$ with $\|\gamma_n-\gamma_{\theta_0,f_0}\|\lesssim \rho_n$ and $n\e_n^4+n\e_n^2\rho_n^2\to0$.
Then the Bernstein-von Mises theorem \eqref{EqBvMTheorem}  holds at $(\q_0,f_0)$.
\end{cor}

\begin{proof}
We apply Theorem~\ref{thm: main} with the sets $\Theta_n=(\q_0-\e_n,\q_0+\e_n)$ and the sets $H_n=H_{n,1}\cap H_{n,2}$, for
$H_{n,1}$ equal to the ball of radius a large multiple of $\e_n$ around $f_0$ and 
$$H_{n,2}=\Bigl\{\forall t\in (-\e_n,\e_n): \Bigl| \log\frac{d\pi_{f+t\gamma_n}} {d\pi_f}(f)\Bigr|\le \eta_n (1+nt^2).\Bigr\}.$$ 
The posterior mass of the sets $\Theta_n$ and $H_{n,1}$ tends to 1 in probability, by assumption.
To see that the same is true for the sets $H_{n,2}$,
we apply the ``remaining mass'' principle (see Theorem~8.20 in \cite{vaartghosal}). The Kullback-Leibler
divergence $K(P_{\q_0,f}^n,P_{\q,f}^n)$ and variance $V_{2,0}(P_{\q_0,f}^n,P_{\q,f}^n)$ are equal to $1/2$ and $1$ times 
$n\|f-f_0\|^2+n\|K_\q f-K_{\q_0}f_0\|^2$, respectively, whence their roots are bounded above by 
$\sqrt n D_3|\q_1-\q_2|\,\|f_0\|+\sqrt n(1+D_3)\|f-f_0\|$,
by Assumption~\ref{asspt:main}(iv). It follows that 
\begin{align*}
&\Pi\bigl((K\vee V_{2,0})(P_{\q_0,f}^n,P_{\q,f}^n) \le n\e_n^2\bigr)\\
&\qquad\ge\pi_\q\Bigl(|\q-\q_0|<\frac{\e_n}{2D_3\|f_0\|}\Bigr) \pi_f\Bigl(\|f-f_0\|<\frac{\e_n}{2D_3+2}\Bigr).
\end{align*}
In view of \eqref{eq: thmPrior} the right side  is bounded below by $e^{-n\e_n^2/65}$, for sufficiently large $n$. 
Combination with \eqref{EqPriorLikelihoodRatio} shows that the remaining mass condition of Theorem~8.20 in \cite{vaartghosal}
is satisfied and hence the posterior probability of the sets $H_{n,2}$ tends to 1.

Thus the first part of \eqref{eq: contraction} is satisfied. For the proof of the second part, we first note that for $\q\in\Theta_n$
the sets $H_{n,1}+(\q-\q_0)\g_n$ are contained in balls of radius a multiple of $\e_n$ around $f_0$ and hence their
posterior probability under $\Pi^{\q=\q_0}$ tends to 1. Secondly, because 
$d\pi_{f+t\g_n}/d\pi_f(f-s\g_n)=d\pi_{f+(s+t)\g_n}/d\pi_{f+s\g_n}(f)$ almost surely, we have 
$$H_{n,2}+s\g_n=\Bigl\{\forall t\in (-\e_n,\e_n): \Bigl| \log\frac{d\pi_{f+(s+t)t\gamma_n}} {d\pi_{f+s\g_n}}(f)\Bigr|\le \eta_n (1+nt^2)\Bigr\}.$$ 
Therefore $\cup_{|s|<\e_n}(H_{n,2}+s\g_n)^c$ is the event in the left side of \eqref{EqPriorLikelihoodRatio}, 
and hence it has prior probability bounded above by $e^{-3n\e_n^2}$. It follows that 
$\sup_{|s|<\e_n}\Pi^{\q=\q_0}\bigl((H_{n,2}+s\g_n)^c | X^n\bigr)
\le \Pi^{\q=\q_0}\bigl(\cup_{|s|<\e_n}(H_{n,2}+s\g_n)^c| X^n\bigr)\ra0$ in probability, by the remaining mass principle.

This verifies the prior shift condition \eqref{eq: priorshift}. The remaining condition \eqref{eq: lfdrate} is first reduced
to \eqref{eq: lfdrate3}, which in turn can be bounded above by $\sqrt n \|f-f_0\|\bigl(\|\g_n-\g_{\q_0,f_0}\|+\|f-f_0\|\bigr)$,
as explained in the discussion preceding the corollary. The latter tends to zero by the rate assumptions $n\e_n^4+n\e_n^2\rho_n^2\to0$.
\end{proof}

As we will see in the next sections, the mapping properties of the solution operator $K_{\theta_0}$ imply that 
$\gamma_{\theta_0, f_0}$ is typically (much) smoother than the true function $f_0$. For this reason it is rarely
useful to apply the preceding corollary with $\rho_n$ slower than $\e_n$. 

On the other hand,  a faster rate $\rho_n$ may be useful, in particular if structural properties of the operators
reduce the ``bias term'' \eqref{eq: lfdrate}. The following corollary drops the requirement
$n\e_n^4\to 0$, under the condition that the range spaces of $K_{\theta_0}$ and its derivative $\dot K_{\theta_0}$ are orthogonal
(or equivalently $\|K_\theta f\| = \|K_\tau f\|$ for all $\theta,\tau\in\Theta$ and $f\in H$), or more generally if \eqref{eq: lfdrate2} holds. 
In this special case, only the condition $\sqrt{n}\eps_n\rho_n\to 0$ on the joint rates remains, in addition to \eqref{eq: priorshift},
and hence the rate of posterior contraction $\e_n$ may be slower than the ``usually'' required rate $o(n^{-1/4})$.

This case includes the case that \eqref{eq: priorshift} is satisfied with the choice $\gamma_n=\gamma_{\theta_0,f_0}$, when
$\rho_n$ can be taken equal to zero and hence no rate of contraction is required.

\begin{cor}\label{cor: constantNorm}
Suppose that either \eqref{eq: lfdrate2} holds or the map $\q\mapsto \|K_\theta f\|$ is constant, for every $f\in H$.
Then the assertion of Corollary~\ref{cor: rateCondition} is true under the stated conditions without the
requirement that $n\e_n^4\to0$.
\end{cor}

\begin{proof}
The proof is the same as the proof of Corollary~\ref{cor: rateCondition}, except that the left side of \eqref{eq: lfdrate}
is bounded by \eqref{eq: lfdrate2}, which vanishes, plus the supremum of 
$\sqrt n \bigl|\ip{f-f_0,(I+K^*_{\theta_0}K_{\theta_0})(\gamma_{\theta_0,f_0}-\gamma_n)}\bigr|$,
which is bounded by $\sqrt n\|f-f_0\|\, \|\gamma_{\theta_0,f_0}-\gamma_n\|\le \sqrt n\e_n\rho_n$.
\end{proof}

A suitable rate of contraction, if needed, can be obtained from standard posterior contraction results, 
as in \cite{vaartghosal_noniid} and specialized to the present model in the following lemma.

\begin{lemma}\label{thm: main2}
Assume that Assumption~\ref{asspt:main} is satisfied and that $\dot K_{\q_0} f_0\not=0$ and $K_\q f_0\not=K_{\q_0}f_0$ 
for all $\q\not=\q_0$.
Suppose there exist $\eps_n\downarrow 0$ with $n\eps_n^2\ra\infty$ 
and  $\FF_n \subset H$ with  $\log \sup_{f\in\FF_n}\|f\|\le  n\eps_n^2/4$ such that \eqref{eq: thmPrior} holds and
\begin{align}
	\pi_f(f\notin\FF_n) &\le e^{-n\eps_n^2},\label{eq: thmPrior2}\\
	\log{N(\eps_n/(8D_3), \FF_n, \|\cdot\|)}&\le  n\eps_n^2/2.\label{eq: thmCover}
\end{align}
Then the posterior distribution of $(\q,f)$ contracts at rate $\e_n$ at $(\q_0,f_0)$, and the same is true for the
posterior  distribution of $f$ given $\q=\q_0$ at $f_0$.
\end{lemma}

\section{Example Applications}
\label{sec:applications}

\subsection{Thermal Diffusivity Recovery in Heat Equation}\label{subsec: heat}

We consider the one-dimensional heat equation which describes the evolution of the temperature in a thin metal 
rod as a function of position and time. The temperature $u(x,t)$ at location $x\in[0,1]$ and $t\in [0,T]$ is governed by the partial differential equation
\begin{align}\label{eq:heatEquation}
    {\frac{\partial }{\partial t}} u(x,t) = \theta \frac{\partial^2}{\partial x^2}u(x,t),\qquad  \text{ \text{ }} u(x,0)=f(x).
    \text{ \text{ }} u(0,t)=u(1,t)=0,
\end{align}
Here $f\in L^2[0,1]$ represents the initial condition of the temperature in the system, and is assumed to satisfy the same boundary conditions $f(0)=f(1)=0$. 
The parameter $\theta > 0$ is the thermal diffusivity constant of the metal and is our main object of interest. 

If $\theta$ is  known, then inferring the initial function $f$ from 
noisy observations of the function $u(\cdot, T)$ at the final time $T$ is a well-known 
statistical inverse problem, see for instance \cite{Bissantz2008}, \cite{Cavalier2008, Cavalier2011}, \cite{Golubev1999}, 
\cite{Knapik2013}, \cite{Mair1994}, \cite{Mair1996}, \cite{Stuart2010}.  The problem is severely ill-posed, in that the solution operator 
$f\mapsto u(\cdot, T)$ takes functions $f\in L^2[0,1]$ into super smooth functions.  
Indeed, if we denote by $f_k$  the  coefficients of $f$ with respect to the sine basis functions  $e_k(x)=\sqrt{2}\sin(k\pi x)$, 
then the solution of  \eqref{eq:heatEquation}  can be expressed  as 
\begin{align}\label{eq: sol}
    u(x,t)=\sum_{k=1}^{\infty}f_k e^{-\theta\pi^2 t k^2}e_k(x).
\end{align}
As a consequence, the optimal, minimax rate for  estimators that recover $f$ from a noisy observation of $u(\cdot,T)$ is known to be of the order $(\log n)^{-\beta/2}$ (for known $\q$ with respect to $L^2$-loss), where $\beta$ is the (Sobolev) regularity of $f$ and $n$ is the signal-to-noise ratio. Various frequentist and Bayesian methods achieve this rate, see for instance \cite{Mair1994, Mair1996, Golubev1999, Bissantz2008, Knapik2013}.

Here, we are interested in learning the diffusivity parameter $\theta$, while still assuming that $f$ is unknown as well.
Just observing the solution of the equation at the final time $T$ is then not sufficient 
to identify the parameter $\theta$, as can be seen from \eqref{eq: sol}. (A change in $\q$ can be fully compensated by a
change in the sequence $f_k$.)  Instead we assume that we have noisy observations of the system 
\eqref{eq:heatEquation} at times $0$ and $T$, i.e.\ of the functions $u(\cdot,0)=f$ and $u(\cdot, T)$.
This is precisely the observational model \eqref{eq: obs1}--\eqref{eq: obs2}, with $H=L^2[0,1]$ and
operator $K_\theta: L^2[0,1]\to L^2[0,1]$ equal to the solution operator 
\begin{equation}\label{eq: k}
K_\theta f = \sum_{k=1}^{\infty}f_k e^{-\theta\pi^2 T k^2}e_k.
\end{equation}
Alternatively, the model observations can be described as a pair of diffusion processes $X^{(i)} = (X^{(1)}_t: t \in [0,T])$, for $i=1,2$,
satisfying, for independent standard Brownian motions $W^{(1)}$ and $W^{(2)}$,
\begin{align}
	\label{eq:sde1} 
dX^{(1)}_t & = f(t)\,dt + \frac1{\sqrt{n}}\,dW^{(1)}_t,\quad
dX^{(2)}_t  = K_\theta f(t)\,dt + \frac1{\sqrt{n}}\,dW^{(2)}_t.
\end{align}
An advantage of the formulation  \eqref{eq: obs1}--\eqref{eq: obs2} is that it extends immediately to the heat equation
on a $d$-dimensional domain $O$, by setting $H=L^2(O)$, letting $e_k$ be the corresponding 
eigenfunctions of the Laplacian and replacing $-(\pi k)^2$ in \eqref{eq: k} by the corresponding eigenvalues.

The map $\theta\mapsto K_{\theta}f$ admits the derivative $\dot K_\theta: L^2[0,1] \to L^2[0,1]$ given by 
\begin{equation}\label{eq: kdot}
\dot K_\theta f =  - \pi^2 T\sum_{k=1}^{\infty}f_k k^2e^{-\theta\pi^2 T k^2}e_k.
\end{equation}
For $\eta\ge 0$, define the Sobolev type space $S^\eta$ with square norm by 
 \begin{equation}
\label{EqDefSeta}
S^\eta:=\Bigl\{f\in H: \sum_{k=1}^\infty f_k^2k^{2\eta}<\infty\Bigr\},
\qquad \|f\|^2_{S^\eta}:= \sum_{k=1}^\infty f_k^2k^{2\eta}.
 \end{equation}

\begin{prop}
Let $H = L^2[0,1]$ and $\Theta = [a,b]\subset (0,\infty)$. Assume that $\eta>1/2$.
Then  the solution operator $K_\theta$ given in \eqref{eq: k} satisfies Assumption~\ref{asspt:main}
relative to $S=S^\eta$. Furthermore \eqref{eq: integralCondition} is satisfied, and $\dot K_{\q_0}f_0\not=0$ 
and $(K_\q-K_{\q_0})f_0\not=0$ provided $f_0\not=0$ and $\q\not=\q_0$.
\end{prop}

\begin{prf}
The entropy of the unit ball of $S^\eta$ in $H$ is known to be of the order $(1/\e)^{1/\eta}$, whence 
\eqref{eq: integralCondition} is satisfied for $\eta>1/2$. The final assertions follow
by injectivity of the operator $\dot K_\q$ and the fact that $(K_\q-K_{\q_0})f=0$ if and only
$(e^{-\q\pi^2 Tk^2}-e^{-\q_0\pi^2 Tk^2})f_k=0$, for all $k$.

For the proof of Assumption~\ref{asspt:main}(i)-(iv) we use the explicit expressions for the operators.
(i). The square norm $\|K_\q f\|^2=\sum_k f_k^2e^{-2\q\pi^2 Tk^2}$ is bounded above by  a multiple of 
$\sum_kf_k^2=\|f\|^2$, in view of the inequality $\sup_{x \ge 0} x^2e^{-x} = 4e^{-2}$. 
(ii). Since $|e^{-x}-e^{-y}| \le |x-y| e^{-x\wedge y}$,
\begin{align*}
\|(K_\theta  - K_{\theta_0})f\|^2_{S^\eta}
 & = \sum_k k^{2\eta}f_k^2(e^{-\theta \pi^2Tk^2}- e^{-\theta_0 \pi^2Tk^2})^2\\
& \le \sum_k k^{2\eta}f_k^2 e^{-2(\theta_0\wedge\theta) \pi^2Tk^2}((\theta-\theta_0) \pi^2Tk^2)^2.
\end{align*}
This is bounded above by a multiple of $\|f\|^2$, since $\sup_{x \ge 0} e^{-x} x^{2+\eta} <\infty$. 
(iii). We have
\begin{align*}
\|(K_{\theta} - K_{\theta_0} -(\theta-\theta_0)\dot K_{\theta_0})f\|^2_{S^\eta} & = 
\sum k^{2\eta}f_k^2 e^{-2\theta_0\pi^2Tk^2}h^2((\theta- \theta_0)\pi^2 T k^2), 
\end{align*}
where $h(x) = e^{-x} -1 + x$. 
The proof is completed as the proof of (ii), using that 
 $|h(x)| \le x^2/2$ for $x \ge 0$ and $\sup_{x \ge 0} e^{-x} x^{4+\eta} = e^{-(4+\eta)}(4+\eta)^{4+\eta}$.
(iv). By the triangle inequality and (i) we have $\|K_{\theta_1} f_1 - K_{\theta_2}f_2\| 
\lesssim \|f_1-f_2\| + \|(K_{\theta_1} - K_{\theta_2})f_2\|$, where the second term is bounded
above by a multiple of $|\q_1-\q_2|\,\|f_2\|$, by a simplified form of the argument under (ii).
Since the same is true with the roles of $f_1$ and $f_2$ reversed, we arrive at the desired inequality. 
\end{prf}

By the explicit expressions for the operators $K_\theta$ and $\dot K_\theta$
the least favorable direction \eqref{EqLeastFavourableDirection} is given by 
\begin{align}\label{eq: lfdHeat}
	\gamma_{\theta_0, f_0} =  - \pi^2 T\sum_{k=1}^{\infty} \frac{k^2 \ip{f_0, e_k}e^{-2\theta_0\pi^2 T k^2}}
	{1+e^{-2\theta_0\pi^2 T k^2}}  e_k.
\end{align}
The exponential decrease of the coefficients of $\gamma_{\theta_0, f_0}$ 
(with respect to the sine basis) implies that this function is analytic.
As a result, as is seen in the next section, Corollary~\ref{cor: rateCondition} can be applied with $\g_n=\g_{\q_0,f_0}$.

\subsection{Location Recovery in Semi-blind Deconvolution}\label{subsec: deconv}
Let $g$ be a known symmetric, square integrable, 1-periodic function with square-integrable derivative.
For $\q\in\Theta\subset (-1/2,1/2)$, consider the convolution operator $K_\theta: L^2[0,1]\to L^2[0,1]$ given by
\[
    K_{\theta}f(t):= g*f(t-\theta)=\int_{0}^{1} f(t-u)g(u-\theta)\,\text{d}u.
\]
Differentiability of $g$ yields, for $g'$ the derivative of $g$,
\[
\dot K_\theta f = -\int_{0}^{1} f(t-u)g'(u-\theta)\,\text{d}u.
\]
An alternative, handy way to write $K_\theta f$ and $\dot K_\theta f$ is through the complex exponential basis $\{e_k(\cdot)\} = \{e^{-2i k\pi\cdot }\}$. Since $g$ is symmetric,
its coefficients $g_k$ in this basis are real and symmetric ($g_k=g_{-k}$). 
Then for $f_k=\ip{f, e_k}$,
\begin{align}\label{eq: kconv}
K_\theta f &= \sum_{k\in\mathbb{Z}}f_k g_k e^{2\pi i k\theta}e_k,\\
\label{eq: kconvDif}
\dot K_\theta f &= -2\pi i\sum_{k\in\mathbb{Z}}f_kk g_k e^{2\pi ik\theta}e_k.
\end{align}
Our objective is to recover $\theta$. As can be seen from \eqref{eq: kconv}, sole knowledge of $K_\theta f$ does not identify $\q$. 
Therefore, we consider noisy observation of $f$ next to $K_\q f$ and consider again model \eqref{eq: obs1}--\eqref{eq: obs2}.
As in the preceding section, the observations can equivalently be defined by the pair of diffusion equations \eqref{eq:sde1}.

Let $S^\eta$ be defined as in \eqref{EqDefSeta}, but with $(f_k)$ the coefficients relative to the present basis
and the index $k$ of the series running through $\mathbb{Z}$.

\begin{prop}\label{prop: deconv}
Let $H=L^2[0,1]$ and let $\Theta\subset (-1/2,1/2)$ be compact.  
Let $\ip{f_0,e_1}\neq 0$, $g_1\neq 0$ and $\sup_kk^{5/2+\d}|g_k|<\infty$, for some $\d>0$. 
Then the solution operator $K_\theta$ given in \eqref{eq: kconv} 
satisfies Assumption~\ref{asspt:main} relative to $S=S^\eta$ for $\eta\le 1/2+\delta$.
Furthermore, \eqref{eq: integralCondition} is satisfied for $\eta>1/2$, the map
$\q\mapsto \|K_\q f\|$ is constant for every $f$,  $\dot K_{\theta_0}f_0\neq 0$ and $(K_\theta-K_{\theta_0})f_0\neq0$ for $\theta\neq \theta_0$.
\end{prop}

\begin{prf}
 (i). By 1-periodicity $\|K_\theta f\|=\|g*f(\cdot-\theta)\| = \|f*g\|$, which is independent of $\q$. Here 
$\|f*g\|\leq \|f\|\,\|g\|$, by the Cauchy-Schwarz inequality (and periodicity), whence 
$\|K_\theta\|=\sup_{\|f\|\le1}\|f*g\|\le \|g\|<\infty$. Similarly $\|\dot K_\theta\|\le \|g'\|<\infty$.
  (ii).  Choose $1/2<\eta \le 1/2+\d$. The sequence $|g_k|k^{\eta+1}$ is thus bounded. Then
 \begin{align*}
 		\|(K_\theta-K_{\theta_0})f\|_{S^\eta}^2 
&= \sum_{k\in\mathbb{Z}}|k|^{2\eta} |f_k|^2g_k^2\,|e^{-i2k\pi\theta}-e^{-i2k\pi\theta_0}|^2\\
&\leq \pi^2|\theta-\theta_0|^2\sum_{k\in\mathbb{Z}}|k|^{2\eta+2}g_k^2|f_k|^2
		\lesssim |\theta-\theta_0|^2\|f\|^2,
 \end{align*}
by the inequality $|e^{ix}-e^{iy}|\le |x-y|$, for every $x,y\in\mathbb{R}$. 
 (iii). As $\eta\leq 1/2+\delta$, the sequence $|g_k|k^{\eta+2}$ is bounded. By the exact formulas
\begin{align*}
 		&\|(K_\theta-K_{\theta_0}-(\theta-\theta_0)\dot K_{\theta_0})f\|_{S^\eta}^2\\
 &\qquad= \sum_{k\in\mathbb{Z}} |k|^{2\eta}|f_k|^2g_k^2\,\bigl|e^{-2\pi i k\theta} - e^{-2\pi i k\theta_0}+(\theta-\theta_0)2\pi k i e^{-i2\pi k\theta_0}\bigr|^2\\
 &\qquad=\sum_{k\in\mathbb{Z}} |k|^{2\eta}|f_k|^2g_k^2\,\bigl|e^{-2\pi i k(\theta-\theta_0)} - 1+2\pi i k(\theta-\theta_0)\bigr|^2.
 \end{align*}
Here $|e^{ix}-1-ix|\leq x^2/2$, for any $x\in\mathbb{R}$, and then  (iii) follows in the same way as (ii).
 (iv). Assume without loss of generality that $\|f_2\| = \|f_1\|\wedge \|f_2\|$. By the triangle inequality,
(i) and (ii) (for $\theta_1$ and $\theta_2$), 
 \begin{align*}
 		\|K_{\theta_1}f_1-K_{\theta_2}f_2\|&\lesssim \|f_1-f_2\|+\|(K_{\theta_1}-K_{\theta_2})f_2\|
\lesssim \|f_1-f_2\| + |\theta_1-\theta_2|\|f_2\|.
\end{align*}

The entropy of the unit ball of $S^\eta$ in $H$ is known to be of the order $(1/\e)^{1/\eta}$, whence 
\eqref{eq: integralCondition} is satisfied for $\eta>1/2$. 

Since both $f_{0,1}:=\ip{f_0,e_1}$ and $g_1$ are nonzero, we have the following non trivial lower bound
  \begin{align*}
  	\begin{split}
  		\|(K_{\theta}-K_{\theta_0})f_0\|^2&\ge g_1 ^2|f_{0,1}|^2|e^{-i2\pi\theta}-e^{-i2\pi\theta_0}|^2\\
  		&= 4g_1 ^2|f_{0,1}|^2\sin\bigl(\pi (\theta-\theta_0)\bigr)^2
		\gtrsim |\theta-\theta_0|^2,
  	\end{split}
  \end{align*}
uniformly in $\q$ such that $\pi|\theta-\theta_0|$ is bounded away from $\pi$. 
This implies that $\dot K_{\q_0}f_0\not=0$ and $K_{\theta}f_0=K_{\theta_0}f_0$ if and only if $\q=\q_0$.
\end{prf}
 
As the convolution operators have constant norm $\q\mapsto \|K_\theta f\|$, they fit in the setting of Corollary~\ref{cor: constantNorm},
meaning that the Bernstein-von Mises theorem may be valid without a rate of posterior contraction.
The least favorable direction \eqref{EqLeastFavourableDirection} has the following expression
\begin{align}\label{eq: lfdDeconv}
	\gamma_{\theta_0, f_0} =  - 2\pi i\sum_{k\in\mathbb{Z}}\frac{k\ip{f_0, e_k}g_k^2}
	{1+g_k^2}  e_k.
\end{align}
Thus  the smoothness of $\gamma_{\theta_0,f_0}$, in the Sobolev sense used in this section,
depends directly on the smoothness of $g$. In the next section we shall see that smoother $g$ 
allow for a wider range of application of the BvM phenomenon. This is interesting, as a smoother $g$ 
means a more ill-posed deconvolution problem, making it harder to recover $f$.


\section{Results for specific priors}
\label{sec:priors}

We present different priors, to which the results of Section~\ref{sec:result} apply, and derive BvM results for the examples of Section~\ref{sec:applications}.

\subsection{Gaussian process priors}
\label{SectionGP}

Let $\pi_f$ be a centered Gaussian process prior in $H$, with reproducing kernel Hilbert space (RKHS) $\HHH$ 
 (see \cite{rkhs}). Define the associated decentering function and concentration functions by, for given $\g, f_0\in H$,
\begin{align*}
\psi_{\g}(\eps) &= \inf_{h \in \HHH:\|h-\g\| \le \eps}\|h\|^2_\HHH,\\
\phi_{f_0}(\eps) &= \psi_{f_0}(\eps)- \log\pi_f\bigl(f: \|f\| \le \eps\bigr).\\
\end{align*}
It is known from \cite{vdVvZ08} that  solutions $\e_n$ to the inequality 
$\phi_{f_0}(\eps_n) \le n\eps^2_n$ give a rate of posterior contraction at $f_0$, in models for which
the intrinsic metric agrees with the norm $\|\cdot\|$. Below we extend this to the present 
model \eqref{eq: obs1}--\eqref{eq: obs2}. It is also known that the prior $\pi_f$ satisfies the
prior mass condition \eqref{eq: thmPrior} (see Theorem~2.1 in \cite{vdVvZ08}) for this value of $\e_n$.

The following application of Corollaries~\ref{cor: rateCondition} and~\ref{cor: constantNorm}
gives BvM results for model \eqref{eq: obs1}-\eqref{eq: obs2} when putting a Gaussian process prior on $f$. 
Cases (ii) and (iv) of the proposition do not involve $\rho_n$. They are obtained from (i) and (iii) by choosing
$\rho_n\asymp n^{-1/2}$, but singled out for special interest.

\begin{prop}\label{prop: GP}
Assume that \eqref{eq: integralCondition} and Assumption~\ref{asspt:main} are satisfied, 
that $\q_0$ is interior to the compact set $\Theta$,  that $\dot K_{\q_0} f_0\not=0$
and that $K_\q f_0\not=K_{\q_0}f_0$ for all $\q\not=\q_0$.
Let $\e_n\downarrow0$ with  $n\eps_n^2\ra\infty$ be such that $\phi_{f_0}(\eps_n)\le n\eps_n^2$, and
let $\rho_n\downarrow0$ be such that $\psi_{\gamma_{\theta_0,f_0}}(\rho_n)\le n\rho_n^2$. 
Suppose that either (i), (ii),  (iii) or (iv) holds:
\begin{enumerate}[label=(\roman*)]
\item $n\eps_n^4+n\e_n^2\rho_n^2\to 0$.
\item $n\eps_n^4\to 0$ and $\g_{\q_0,f_0}\in\HHH$.
\item $\q\mapsto \|K_\theta f\|$ is constant on $\Theta$, and $n\eps_n^2\rho_n^2\to 0$.
\item $\q\mapsto \|K_\theta f\|$ is constant on $\Theta$, and $\g_{\q_0,f_0}\in\HHH$.
\end{enumerate}
Then the Bernstein-von Mises theorem \eqref{EqBvMTheorem} holds at $(\q_0,f_0)$.
\end{prop}

\begin{proof}
We verify the conditions of Corollary~\ref{cor: rateCondition} for (i), and of Corollary~\ref{cor: constantNorm} for (iii).
Cases (ii) and (iv) follow from (i) and (iii) by choosing $\rho_n\asymp n^{-1/2}$.

Condition \eqref{eq: integralCondition} and Assumption~\ref{asspt:main} are satisfied, and  $\dot K_{\q_0} f_0\not=0$, by assumption.
Posterior contraction at the rate $\e_n$ follows by Lemma~\ref{thm: main2}. 
Condition \eqref{eq: thmPrior} is satisfied (for a multiple of $\e_n$), by Theorem~2.1 in \cite{vdVvZ08}.
The rate conditions on $\e_n$ or $\rho_n$ (if any) are valid by assumption.
It remains to verify \eqref{EqPriorLikelihoodRatio}.

For any $\g_n\in\HHH$, the shifted prior $\pi_{f+t\g_n}$ is absolutely continuous relative to $\pi_f$ with density
$$\frac{d \pi_{f+t\g_n}}{d\pi_f}(f)=e^{t \|\g_n\|_\HHH U(f)-t^2\|\g_n\|_\HHH^2},$$
for a measurable transformation $U(f)$ with a standard normal distribution, if $f\sim\pi_f$ (see Lemma~3.2 in \cite{rkhs}). 
By the definition of the decentering function and the assumption on $\rho_n$, there exist
$\gamma_n\in \HHH$ such that $\|\gamma_n-\gamma_{\theta_0,f_0}\|\lesssim \rho_n$ and $\|\gamma_n\|^2_\HHH\leq n\rho_n^2$.
By applying the preceding display twice,
\begin{align*}
\frac{\bigl|\log (d \pi_{f+(s+t)t\g_n}/d\pi_{f+s\g_n})(f)\bigr|}{1+nt^2}&
\le\frac{\|\g_n\|_\HHH|U(f;\g_n)|}{\sqrt n}+\frac{2|s|\|\g_n\|_\HHH^2}{\sqrt n}+\frac{\|\g_n\|_\HHH^2}{n}\\
&\lesssim \|\g_n\|_\HHH\e_n+\frac{\e_n\|\g_n\|_\HHH^2}{\sqrt n}+\frac{\|\g_n\|_\HHH^2}{n},
\end{align*}
on the event $\{f: |U(f)|\ge 2\sqrt n\e_n\}$, for $|s|<\e_n$. By the tail bound on the standard normal distribution, the complement
of the latter event has $\pi_f$-probability smaller than $e^{-3n\e_n^2}$. This verifies \eqref{EqPriorLikelihoodRatio}
provided the right side of the preceding display tends to 0, i.e.\ $\sqrt n \rho_n\e_n\ra0$ and $\rho_n\ra 0$.

If $\g_{\q_0,f_0}\in\HHH$, then we choose $\g_n=\g_{\q_0,f_0}$ and \eqref{EqPriorLikelihoodRatio} is
easily verified by the preceding argument, for any $\e_n\ra0$.
\end{proof}

We now consider natural and commonly used choices of centered Gaussian process priors 
$\pi_f$. These priors all have a hyper-parameter that can be viewed as describing a form of 
``smoothness'', or ``regularity'' of the prior. 
We investigate in particular for which combinations of prior regularity and regularity of the 
true function $f_0$ we can apply Proposition~\ref{prop: GP}(i) to the heat equation example (Section~\ref{subsec: heat})
and Proposition~\ref{prop: GP}(iii) to the deconvolution example (Section~\ref{subsec: deconv}).

\begin{ex}[Series prior]\label{ex: seriesprior}
	Since the solution operator  of the heat equation \eqref{eq: k} diagonalizes on the sine basis, and the convolution operator \eqref{eq: kconv} diagonalizes on the complex exponential basis, it is natural to 
	consider priors on $f$ with covariances that diagonalize on those bases as well. 
	In this example we therefore consider the prior  $\pi_f$ 
	defined as the distribution of the random series 
	\begin{equation*}
		f = \sum_k \sigma_kZ_ke_k,
	\end{equation*}
	where the $e_k$ are either the sine basis functions (for the heat equation example) or the complex exponentials (for the deconvolution example), the $Z_k$ are independent standard normal variables and $\sigma_k$ is a sequence of standard deviations that satisfies $\sum_k \sigma^2_k < \infty$, 
	ensuring that the random series defines a random element in $L^2[0,1]$. In particular we consider $\sigma_k \asymp |k|^{-1/2-\alpha}$, for $\alpha> 0$.  
	This yields a prior on $f$  which   (almost) has regularity $\alpha$ in the Sobolev-type sense with respect 
	to the basis $\{e_k\}$:  $\sum k^{2s}\qv{f, e_k}^2 < \infty$,	almost surely, for all $s < \alpha$.
	We assume that $f_0$ has regularity $\beta$ in the same sense, that is, $\|f_0\|^2_{S^\beta} = 
	\sum |k|^{2\beta}|\qv{f_0, e_k}|^2 < \infty$. 

The reproducing kernel Hilbert space of the process $f$ given by (see Theorem~4.2 of \cite{rkhs})
	\begin{equation*}
		\HHH = \Big\{h= \sum_k c_k e_k: \|h\|^2_\HHH = \sum_k \frac{|c_k|^2}{\sigma^2_k} < \infty\Big\}. 
	\end{equation*}
The function $h_K = \sum_{|k| \le K}\qv{f_0, e_k}e_k $ is contained in $\HHH$, with 
$\|h_K - f_0\| \le K^{-\beta}\|f_0\|_{S^\beta}$ and 	$\|h_K\|^2_\HHH = \sum_{|k| \le K} |\qv{f_0, e_k}|^2/\sigma^2_k \lle 
	(1\vee K^{1+2(\alpha - \beta)})\|f_0\|^2_{S^\beta}$. Therefore, given $\eps > 0$,
 we can take $K \asymp (\|f_0\|_{S^\beta}/\eps)^{1/\beta}$ in the definition of the decentering function, to find that 
$\psi_{f_0}(\eps)  \lle \eps^{-{1+2(\alpha - \beta)}/{\beta}} \vee 1$.
Combining this with Corollary~4.3 of \cite{dunker}, we find that $\phi_{f_0}(\eps)=\psi_{f_0}(\eps)  -\log\pi_f(\|f\| \le \eps) \lle  
\eps^{-{1+2(\alpha - \beta)}/{\beta}} \vee 1+ \eps^{-1/\alpha}$.
Hence  it follows that the inequality $\phi_{f_0}(\eps_n) \le n\eps^2_n$ is solved for 
$\eps_n \asymp n^{-(\alpha \wedge \beta)/(1+2\alpha)}$. The condition $n\eps^4_n \to 0$ translates into the 
	requirement that $1/2 < \alpha < 2\beta - 1/2$. 
	
	In the heat equation example the coefficients of $\gamma_{\theta_0, f_0}$ in \eqref{eq: lfdHeat}
	with respect to the sine basis decay exponentially fast. This implies that $\gamma_{\theta_0, f_0} \in \HHH$, irrespective
	of the values of $\alpha$ and $\beta$. Hence Proposition~\ref{prop: GP}(ii) gives the Bernstein-von Mises result
 if $1/2 < \alpha < 2\beta - 1/2$.

In the deconvolution example, the coefficients of $\gamma_{\theta_0,f_0}$ in \eqref{eq: lfdDeconv} 
with respect to the exponential basis are of order $g_k^2\cdot|\qv{f_0,e_k}|\cdot |k|$. If $g_k$ is of exponential order (e.g. if $g$ is a Gaussian kernel), then $\gamma_{\theta_0,f_0}\in \HHH$ no matter the values of $\alpha$ and $\beta$, 
so we are in scenario (iv) of Proposition~\ref{prop: GP} and obtain the BvM result without any rate condition. 
If $g_k$ is of polynomial order $|k|^{-p}$, then  
$\gamma_{\theta_0,f_0}\in S^{\beta+2p-1}$ and we require $p\ge3$ in order to satisfy the conditions of
Proposition~\ref{prop: deconv}. Similarly as above, we deduce that 
$\psi_{\gamma_{\theta_0,f_0}}(\rho_n)\le n\rho_n^2$ for $\rho_n\asymp n^{-(\beta+2p-1)/(1+2\alpha)}$
(cf.\ Lemma~11.41 in \cite{vaartghosal}). 
The BvM result of Proposition~\ref{prop: GP} holds if $n\e_n^2\rho_n^2\ra0$, i.e.
$2(\a\wedge \b)+2(\b+2p-1)>1+2\a$. This translates into $\alpha < 2\beta + 2p -3/2$.
\end{ex}
	
	The next examples show that it is not necessary to use a prior on $f$ that is compatible
	with the operator $K_\theta$. Commonly used Gaussian process priors work 
	just as well and give similar results. Only the appropriate type of regularity to describe 
	the result is slightly different in each case. 
	
\begin{ex}[Integrated Brownian motion]
	Define the $k$-fold integration operators $I_{0+}^k$ by $I_{0+}^0 f = f$ and  $I^{k}_{0+}f(t) 
	= \int_0^t I^{k-1}_{0+}f(s)\,ds$, for $k \in \NN$.  Now fix $k \in \NN_0$ and let $\pi_f$ be the law 
	of the centered Gaussian process 
	\[
	f(x) = \sum_{i=0}^{k} x^iZ_i + I_{0+}^k B(x), \qquad x \in [0,1],  
	\] 
	where  $Z_0, \ldots, Z_{k}$ are independent standard normal variables and $B$ is a standard Brownian motion independent of the $Z_i$. 
	(The independent polynomial is added to the $k$-fold integrated Brownian motion  because the process itself
	and its $k$ derivatives would otherwise all vanish at $0$, which is undesirable.) 
	The well-known properties of Brownian motion imply that the process $W$ (almost) has regularity 
	$\alpha = k+1/2$, in the sense that its sample paths  belong to $C^s[0,1]$, almost surely, for every $s < \alpha$. 
	We assume that $f_0 \in C^\beta[0,1]$ for $\beta > 0$.
	
	It is known that in this situation the concentration inequality $\phi_{f_0}(\eps) \le n\eps^2_n$
	is satisfied for $\eps_n \asymp n^{-\alpha \wedge \beta/(1+2\alpha)}$, see Section 11.4.1 of \cite{vaartghosal}. 
	The RKHS of $W$ is the (usual) $L^2$-Sobolev space of regularity $\alpha+1/2 = k+1$, consisting  
	of the functions on $[0,1]$ that are $k$ times differentiable with $k$th derivative $f^{(k)}$ 
	that is absolutely continuous  with derivative $f^{(k+1)}$ in $L^2[0,1]$ (see Lemma~11.29 of
	\cite{vaartghosal}). 
	
In the heat equation example $\gamma_{\theta_0, f_0}$ in \eqref{eq: lfdHeat} is infinitely often continuously differentiable.
Therefore, it belongs to the RKHS and Proposition~\ref{prop: GP} (ii) shows that the Bernstein-von Mises theorem holds,
under only the condition $n\e_n^4\to0$, i.e.\ for  if $1/2 < \alpha < 2\beta - 1/2$.
	
In the deconvolution example, differentiability of $\gamma_{\theta_0, f_0}$ in \eqref{eq: lfdDeconv} depends directly on the convolution kernel $g$. If the Fourier coefficients $g_k$ decrease exponentially fast, then $\gamma_{\theta_0,f_0}$ is infinitely often continuously differentiable and the BvM result of Proposition~\ref{prop: GP} (iv) holds without any restriction on $\alpha$ or $\beta$. If this is not the case, then we assume that $g\in C^p[0,1]$. We assume that $p\geq 3$, as this is more than enough to satisfy the conditions of Proposition \ref{prop: deconv}. Observe further that we can rewrite  $\gamma_{\theta_0, f_0} = (I+K_{0}^2)^{-1}g*g'*f_0$. A straightforward analytical argument on convolutions of Hölder continuous functions then shows that $g*g'*f_0\in C^{\beta+2p-1}[0,1]$ and then also
that $\gamma_{\theta_0, f_0}\in C^{\beta+2p-1}[0,1]$, by the mapping properties of the (Fredholm) operator $I+K_{0}^2$.
It follows that $\psi_{\gamma_{\theta_0,f_0}}(\rho_n)\le n\rho_n^2$ for $\rho_n\asymp n^{- (2p-1+\beta)/(1+2\alpha)}$, and the BvM result of Proposition~\ref{prop: GP} (iii) also holds in this case if $1/2 < \alpha < 2\beta + 2p -3/2.$
	
	It is straightforward to extend this example from multiply integrated Brownian motion to the more 
	general Riemann-Liouville process, which also covers fractional integrals of Brownian motion. See 
	Section 11.4.2 of \cite{vaartghosal}.
\end{ex}

\begin{ex}[Mat\'ern process]
	In this example we let $\pi_f$ be the law of a (one-dimensional) Mat\'ern process with parameter 
	$\alpha > 0$. This is a  centered, stationary Gaussian process  $W$ with spectral measure 
	$\mu_\alpha(d\lambda) = (1+\lambda^2)^{-1/2-\alpha}d\lambda$, that is, 
	\[
	\EE W_s W_t = \int_\RR \frac{e^{i\lambda(t-s)}}{(1+\lambda^2)^{1/2+\alpha}}\,d\lambda. 
	\]
	(The stationary Ornstein-Uhlenbeck process is a particular example, corresponding to $\alpha = 1/2$.)
	It can  be shown that the sample paths of the Mat\'ern process almost surely belong to $C^s[0,1]$ 
	for all $s < \alpha$. We assume that for $\beta > 0$ we have $f_0 \in C^\beta[0,1] \cap H^\beta[0,1]$, 
	where $H^\beta[0,1]$ is  the Sobolev space defined as the space of restrictions to $[0,1]$ 
	of functions $f \in L^2(\RR)$ with Fourier transform $\hat f$ that satisfies $\int (1+\lambda^2)^\beta
	|\hat f(\lambda)|^2\,d\lambda < \infty$. 
	
	By Lemmas 11.36 and 11.37 of \cite{vaartghosal}, we have $\phi_{f_0}(\eps_n) \le n\eps^2_n$
	for $\eps_n \asymp cn^{-\alpha \wedge \beta/(1+2\alpha)}$.
	The RKHS of $W$ is the space of (restrictions to $[0,1]$ of) real parts of functions $h$ that can be written 
	as $h(t) = \int e^{i\lambda t}\psi(\lambda)\,\mu_\alpha(d\lambda)$ for some $\psi \in L^2(\mu_\alpha)$, 
with RKHS norm given by $\|h\|_\HHH = \|\psi\|_{L^2(\mu)}$ (see Lemma 11.35 of \cite{vaartghosal}).

	In the heat equation example, smoothness on $[0,1]$ implies that $\gamma_{\theta_0, f_0}$ can be extended to a compactly supported 
	$C^\infty$-function on the whole line. The Fourier transform $\hat \gamma_{\theta_0, f_0}$ of 
	this extension then has the property that $|\lambda|^p |\hat \gamma_{\theta_0, f_0}(\lambda)| \to 0$
	for every $p > 0$ as $|\lambda| \to \infty$. By Fourier inversion the extended function
	can be written as 
	\[
	\gamma_{\theta_0, f_0}(t) = \int e^{i\lambda t} \hat \gamma_{\theta_0, f_0}(\lambda)\,d\lambda 
	= \int e^{i\lambda t}\psi(\lambda)\mu_\alpha(d\lambda), 
	\]
	where $\psi(\lambda) = (1+\lambda^2)^{1/2+\alpha}\hat \gamma_{\theta_0, f_0}(\lambda)$. By the observation 
	just made about the tails of $\hat \gamma_{\theta_0, f_0}$, we have that $\psi \in L^2(\mu_\alpha)$, whence $\gamma_{\theta_0, f_0}$ belongs to the RKHS 
	of $W$. As in the preceding examples the condition on $\gamma_{\theta_0, f_0}$ is trivially 
	satisfied and the Bernstein-von Mises result holds if $1/2 < \alpha < 2\beta - 1/2$,  by Proposition~\ref{prop: GP} (ii).
	
	In the deconvolution example, if the Fourier transform $\hat g$ of the convolution kernel $g$ has the property that $|\lambda|^p|\hat g(\lambda)|\to 0$ for any $p>0$, then the same holds for the Fourier transform $\hat \gamma_{\theta_0,f_0}$ of $\gamma_{\theta_0,f_0}$. It follows that $\gamma_{\theta_0, f_0}$ belongs to the RKHS of $W$ and the BvM result holds without any restrictions on $\alpha$ and $\beta$, by Proposition~\ref{prop: GP} (iv). 
If this is not the case, we assume that $g\in C^p[0,1] \cap H^p[0,1]$ for some $p$. We again have that  $\psi_{\gamma_{\theta_0,f_0}}(\rho_n)\le n\rho_n^2$ for $\rho_n\asymp n^{-(2p-1+\beta)/(1+2\alpha)}$. As for the other priors, the BvM result holds once again in this case  $1/2 < \alpha < 2\beta + 2p -3/2$,
in view of  of Proposition~\ref{prop: GP} (iii).
	
\end{ex}

\begin{ex}[Squared exponential process]
	Let $\pi_f$ be the law of the squared exponential process $W$ with length scale
	$l = l_{n, \alpha} =  n^{-1/(1+2\alpha)}$ for $\alpha > 0$, so 
	\[
	\EE W_s W_t = e^{-(t-s)^2/l^2}.
	\]  
	The sample paths of the squared exponential process are analytic functions, but 
	still in view of the results of \cite{scaling} it makes sense to think of the rescaled process 
	$W$ as ``essentially'' having regularity $\alpha$. We assume that $f_0 \in C^\beta[0,1]$
	for $\beta > 0$. 
	
	By Lemma 2.2 and Theorem 2.4 of \cite{scaling} the inequality $\phi_{f_0}(\eps_n) \le n\eps^2_n$
	holds in this case  for 
	$\eps_n \asymp n^{-\alpha \wedge \beta/(1+2\alpha)}\log n$..
	Also in this case we have $n\eps_n^4 \to 0$ if and only if 
	$1/2 < \alpha < 2\beta - 1/2$, the extra logarithmic factor in the rate $\eps_n$ has no influence 
	on this condition.

	The spectral measure of the stationary process $W$ is given by $\mu(d\lambda) = \mu_{n,\alpha}(d\lambda) 
	= l(2\sqrt\pi)^{-1}e^{-l^2\lambda^2/4}\,d\lambda$
	and as in the preceding example the RKHS $\HHH$ of $W$ is the space of (restrictions to $[0,1]$ of) real parts of functions $h$ that can be written 
	as $h(t) = \int e^{i\lambda t}\psi(\lambda)\mu(d\lambda)$ for some $\psi \in L^2(\mu)$, 
	with RKHS norm given by $\|h\|_\HHH = \|\psi\|_{L^2(\mu)}$
	(see Lemma 11.35 of \cite{vaartghosal}). 
	
	In the preceding example, we noted that $\gamma_{\theta_0, f_0}$ coming from the heat equation example \eqref{eq: lfdHeat} extends to a compactly support 
	$C^\infty$-function on $\RR$ that can be written as 
	$\gamma_{\theta_0, f_0}(t) = \int e^{i\lambda t} \hat \gamma_{\theta_0, f_0}(\lambda)\,d\lambda$,
	where the Fourier transform $\hat \gamma_{\theta_0, f_0}(\lambda)$ has tails that decay faster 
	than any polynomial. Now for $K > 0$, let 
	\[
	h(t) = \int_{-K}^K e^{i\lambda t} \hat \gamma_{\theta_0, f_0}(\lambda)\,d\lambda.
	\]
	Then for every $p> 0$ we have   
	$\|h-\gamma_{\theta_0, f_0}\|^2 
	\lle \int 1_{|\lambda| > K}|\lambda|^{-1-2p}\,d\lambda \lle K^{-2p}$.
	Moreover we have $h(t) = \int e^{i\lambda t}\psi(\lambda)\mu(d\lambda)$, where
	\[
	\psi(\lambda) = \frac{2\sqrt\pi}{l}1_{|\lambda| \le K} e^{\frac{l^2\lambda^2}{4}}\hat \gamma_{\theta_0, f_0}(\lambda).
	\]
	It follows that $\|h\|^2_\HHH = \|\psi\|^2_{L^2(\mu)} \lle l^{-1} e^{\frac{l^2K^2}{4}}$.
	Taking  $K = 1/l$ and $p > \alpha \wedge \beta$, we obtain
	\[
	\inf_{h \in \HHH:\|h-\gamma_{\theta_0, f_0}\| \le \eps_n}\|h\|^2_\HHH \lle n^{1/(1+2\alpha)} \le n\eps^2_n
	\]
	and hence, after enlarging the constant $c$ if necessary,  
	$\psi_{\gamma_{\theta_0, f_0}}(\eps_n) \le n\eps^2_n$. 
	So also with this prior $\pi_f$ the Bernstein-von Mises result of holds if $1/2 < \alpha < 2\beta - 1/2$, by Proposition~\ref{prop: GP} (i).
	
	In the deconvolution example, similar computations show that if the tails of the Fourier transform of the convolution kernel $g$ decay exponentially fast, then $\psi_{\gamma_{\theta_0,f_0}}(\rho_n)\le n\rho_n^2$ for $\rho_n$ of the order $n^{-\alpha/(1+2\alpha)}\log n$. We obtain the BvM result of Proposition~\ref{prop: GP} (iii) if and only if $\alpha\wedge\beta > 1/2$. Otherwise, we assume that $g\in C^p[0,1]$ and get that the inequality holds for $\rho_n \asymp n^{-\alpha\wedge (2p-1+\beta)/1+2\alpha}\log n$. The logarithmic factor does not influence the restrictions on $\alpha, \beta, p$. We obtain that the BvM result of Proposition~\ref{prop: GP} holds for  $\beta>1/2$ and $1/2 < \alpha < 2\beta + 2p -3/2.$
	
	By the results of \cite{scaling} this example generalizes to any 
	process $W_{t/l}$, where $W$ is a centered stationary Gaussian process  with a spectral 
	measure $\mu$ that satisfies $\int e^{\delta|\lambda|}\,\mu(d\lambda) < \infty$ for some $\delta > 0$. 
\end{ex}

The preceding examples show that many common choices of Gaussian process priors $\pi_f$ lead to  
the Bernstein-von Mises result for the marginal posterior of $\theta$. 

Although the appropriate notion of smoothness depends on the prior,
in the example of the heat equation the regularity $\alpha$ of the prior 
and the regularity $\beta$ of the true $f_0$ should satisfy the constraints 
$1/2 < \alpha < 2\beta - 1/2$. Thus the prior $\pi_f$ need not be tuned 
for obtaining an optimal contraction rate for $f$. Both prior undersmoothing, and a limited degree
of oversmoothing are permitted. This phenomenon has been observed also
in other instances of the semi-parametric Bernstein-von Mises theorem, see 
for instance \cite{castillo}, \cite{rene}, \cite{ismaeljudith}. 

In the deconvolution example the permitted range is $1/2<\alpha<2\beta+2p-3/2$.
This shows that the ``smoothening" effect of the operator $K_{\theta_0}$ 
has a direct effect on the Bernstein von-Mises phenomenon. The smoother the convolution kernel, 
the more the prior can oversmooth the truth, as the upper bound $2\beta+2p-3/2$
 increases with $p$. In the case of a super smooth kernel, this restriction disappears completely 
and the BvM holds for virtually any Gaussian prior on $f$.  

\subsection{$p$-exponential priors}
Let  $e_1,e_2,\ldots$ be a given orthonormal basis for $H$ and let $\sigma_k$ be positive numbers with  $\sum_{k=1}^\infty\sigma_k^2<\infty$. 
For $p\in [1,2]$ the $p$-exponential prior $\pi_f$ on $H$  is defined as the distribution of the random series
\[
f = \sum_{k=1}^\infty \sigma_k Z_k e_k,
\]
where the $Z_k$ are i.i.d.\ real random variables with density proportional to $z\mapsto e^{-{|z|^p}/{p}}$.
 The $2$-exponential prior
is a Gaussian process prior, as in Section~\ref{SectionGP}. Theory for general $p$-exponential priors is
developed in \cite{pexp}. In the following we choose $\sigma_k=k^{-1/2-\alpha/d}$, for $\alpha>0$ and $d\in\mathbb{N}$, so as to make
easy comparison to $d$-dimensional Gaussian process priors of smoothness $\alpha$ in $d$ dimensions.

Let $h_k=\ip{h,e_k}$ be the coefficients of $h\in H$ relative to the basis $e_k$.
Associated to $\pi_f$ are the two weighted sequence spaces
\begin{align*}
\mathcal{Z}&=\Bigl\{h\in H: \sum_{k=1}^\infty \frac{|h_k|^p}{\sigma_k^{p}} <\infty \Bigr\},\qquad
\|h\|_{\mathcal{Z}}:=\Bigl(\sum_{k=1}^\infty \frac{|h_k|^p}{\sigma_k^{p}} \Bigr)^{1/p},\\
\mathcal{Q}&=\Bigl\{h\in H: \sum_{k=1}^\infty \frac{|h_k|^2}{\sigma_k^{2}} <\infty \Bigr\}, \qquad
\|h\|_{\mathcal{Q}}:=\Bigl(\sum_{k=1}^\infty \frac{|h_k|^2}{\sigma_k^{2}} \Bigr)^{1/2}.
\end{align*}
For $p=2$ the two spaces coincide and are equal to the reproducing kernel Hilbert space of the Gaussian prior. 


We define decentering and concentration functions associated to  $\pi_f$ by, for $\g,f_0\in H$,
\begin{align*}
\psi_{\g}(\eps) &= \inf_{h \in \mathcal{Z }:\|h-\g\| \le \eps}\|h\|^p_{\mathcal{Z}},\\
\phi_{f_0}(\eps) &= \inf_{h \in \mathcal{Z}:\|h-f_0\| \le \eps}\frac{1}{p}\|h\|^p_{\mathcal{Z}} - \log\pi_f\bigl(f: \|f\| \le \eps\bigr).
\end{align*}
It is shown in \cite{pexp} that solutions $\eps_n$ to $\phi_{f_0}(\eps_n) \le n\eps^2_n$ such that $n\eps^2_n \to \infty$ 
give a rate of posterior contraction in standard models. It can be shown that
this includes the present model \eqref{eq: obs1}-\eqref{eq: obs2}.
In the case that $\sigma_k=k^{-1/2-\alpha/d}$ and $f_0\in S^{\b/d}$, for $S^\eta$ defined in \eqref{EqDefSeta}, 
Proposition~5.4 in \cite{pexp} shows that the solution  to $\phi_{f_0}(\eps_n) \le n\eps^2_n$ is given by 
\begin{equation}
\label{EqRatepExp}
\e_n:=r_n^{\alpha,\beta, p}\asymp\begin{cases} n^{-\beta/(d+2\beta+p(\alpha-\beta))}&\text{ if }\alpha\ge \beta, \\ 
n^{-\alpha/(d+2\alpha)} &\text{ if }\alpha <\beta.\end{cases}
\end{equation}
It is also shown in  \cite{pexp} that \eqref{eq: thmPrior} holds for (a multiple of) the same $\e_n$.

Thus for an application of Theorem~\ref{thm: main} or Corollaries~\ref{cor: rateCondition} or~\ref{cor: constantNorm}, it suffices
to verify the prior shift condition \eqref{eq: priorshift}, or \eqref{EqPriorLikelihoodRatio}. It is known that
$\pi_{f+\g}$ is absolutely continuous with respect to $\pi_f$ if and only if the shift $\g$ belongs to $\mathcal{Z}$.
If the least favorable direction $\gamma_{\theta_0,f_0}$ is not included in $\mathcal{Z}$, then we approximate it by 
a sequence $\g_n\in \mathcal{Z}$ with some approximation rate $\|\gamma_n-\gamma_{\theta_0,f_0}\|\lesssim \rho_n$. 
The decentering function $\psi_{\gamma_{\theta_0,f_0}}(\rho_n)$ 
then controls the rate of growth of $\|\g_n\|_{\mathbb{Z}}$, which allows to verify \eqref{eq: priorshift},
as shown by the following lemma. The proof of the lemma is given in Section~\ref{sec: proof}.

\begin{lemma}
\label{lemma: priorshift}
Let $\e_n$ be positive numbers with $n\e_n^2\gtrsim 1$ so that $\phi_{f_0}(\eps_n) \le n\eps^2_n$ and let $(\gamma_n)\subset\mathcal{Z}$. 
Then \eqref{eq: priorshift} is verified for sets $\Theta_n\times H_n$ satisfying \eqref{eq: contraction} if 
\begin{enumerate}[label=(\roman*)]
\item $p=1$ and $n^{-1/2}\|\gamma_n\|_{\mathcal{Z}}\to 0$.
\item $p\in (1,2]$ and $\eps_n\|\gamma_n\|_{\mathcal{Q}} +n^{-1/2}\eps_n^{p-1}\|\gamma_n\|_{\mathcal{Z}}^p+n^{-1/2}\|\gamma_n\|_{\mathcal{Z}}\to 0$.
\end{enumerate}
A sufficient condition in case (ii), $p\in (1,2]$, is that $\eps_n\|\gamma_n\|_{\mathcal{Z}} \to 0$.
\end{lemma}

If the least favorable direction $\g_{\q_0,f_0}$ is contained in $\mathcal{Z}$, then we can choose $\g_n=\g_{\q_0,f_0}$ and the conditions 
in the lemma are trivially satisfied. Then the only remaining condition in Corollary~\ref{cor: rateCondition} is $n\e_n^4\to0$, while in Corollary~\ref{cor: constantNorm}
also this condition can be dropped. 

If the least favorable condition is not contained in $\mathcal{Z}$, then the conditions are more involved.
By the definition of the decentering function there exist $\g_n$ with $\|\g_n-\g_{\q_0,f_0}\|\le\rho_n$ 
such that $\|\gamma_n\|_{\mathcal{Z}}^p \le \psi_{\g_{\q_0,f_0}}(\rho_n)$.
Thus in the case that $p=1$, it suffices that $n^{-1/2}\psi_{\g_{\q_0,f_0}}(\rho_n)\to0$, next to the condition that $n\e_n^2\rho_n^2\ra0$.
The case $p\in(1,2)$, as in (ii) of the lemma, involves two different norms of $\g_n$, but the sufficient condition $\eps_n\|\gamma_n\|_{\mathcal{Z}} \to 0$
can be captured in the modulus as $\e_n^2 \psi_{\g_{\q_0,f_0}}(\rho_n)^{2/p}\to0$ and together with the condition
 $n\e_n^2\rho_n^2\ra0$ leads to optimizing $\e_n^2\bigl[n\rho_n^2+\psi_{\g_{\q_0,f_0}}(\rho_n)^{2/p}\bigr]$ over $\rho_n$, as in (iii) of the
following proposition.

\begin{prop}\label{prop: pexp}
Assume that \eqref{eq: integralCondition} and Assumption~\ref{asspt:main} are satisfied, 
that $\q_0$ is interior to the compact set $\Theta$,  that $\dot K_{\q_0} f_0\not=0$
and that $K_\q f_0\not=K_{\q_0}f_0$ for all $\q\not=\q_0$.
Let $\pi_f$ be a $p$-exponential prior with scaling sequence $\sigma_k\asymp k^{-1/2-\alpha/d}$. 
Let $\e_n$ be positive numbers with $\e_n\ra0$ and $n\e_n^2\ra\infty$ such that $\phi_{f_0}(\eps_n)\le n\e_n^2$.
Assume one of the following:
\begin{enumerate}[label=(\roman*)]
\item $\gamma_{\theta_0,f_0}\in\mathcal{Z}$ and $n\eps_n^4\to 0$.
\item $p=1$ and $n\eps_n^4\to 0$ and there exists $\rho_n\ra0$ such that $n\eps_n^2\rho_n^2+n^{-1/2}\psi_{\gamma_{\theta_0,f_0}}(\rho_n)\to0$.
\item $p>1$ and $n\e_n^4\to0$ and $n\eps_n^2\rho_n^2\to0$ for $\rho_n$ satisfying $\psi_{\gamma_{\theta_0,f_0}}(\rho_n)\le n^{p/2}\rho_n^p$.
\end{enumerate}
Then  the Bernstein-von Mises theorem \eqref{EqBvMTheorem} holds at $(\q_0,f_0)$.
If the map $\q\mapsto\|K_\q f\|$ is constant on $\Theta$, for every $f\in H$, then (i)--(iii) remain true without the condition that $n\e_n^4\to0$.
\end{prop}

\begin{ex}[Heat equation]
In the example of the heat equation in Section~\ref{subsec: heat}, the Fourier coefficients of 
the least favorable direction $\gamma_{\theta_0, f_0}$ decay exponentially, as seen in \eqref{eq: lfdHeat}.
This implies that $\gamma_{\theta_0,f_0}\in \mathcal{Z}$, if $\{e_k\}$ is chosen equal to the sine basis. For $f_0$ of regularity
$\beta>0$, the posterior rate of contraction is given by \eqref{EqRatepExp}. The condition $n\eps_n^4\to 0$ then
translates into $1/2<\alpha<\frac{2+p}{p}\beta - \frac{1}{p}$. The BvM is valid for $\alpha$  in this range, by
Proposition~\ref{prop: pexp} (i). We observe in particular that a smaller $p$ means a broader application of the BvM.
\end{ex}

\begin{ex}[Convolution]
In the deconvolution problem of Section~\ref{subsec: deconv}, the operator satisfies the condition that $\q\mapsto\|K_\q f\|$ is constant.
Thus the condition that $n\e_n^4\to0$ may be dropped from Proposition~\ref{prop: pexp}, in view of its final assertion.

If the convolution kernel $g$ has exponentially decaying Fourier coefficients with respect to the complex exponential basis, 
then so does the least favorable direction $\gamma_{\theta_0,f_0}$, given in \eqref{eq: lfdDeconv}. 
So $\gamma_{\theta_0,f_0}\in \mathcal{Z}$, for $e_k$ the complex exponential basis in this case, and we get the BvM result 
of Proposition~\ref{prop: GP} (i) without any rate condition. 

If the Fourier coefficients of $g$ are of polynomial order $|k|^{-q}$, then we require $q>5/2$ for Proposition~\ref{prop: deconv} to be applicable. 
(We changed notation, denoting the order of the Fourier coefficients now by $q$ instead of $p$, as $p$ refers to the prior here.) 
We assume that $f_0\in S^\beta$ and then have that $\gamma_{\theta_0,f_0}\in S^{\beta+2q-1}$, in view of \eqref{eq: lfdDeconv}. 
Lemma 5.13 in \cite{pexp} then gives 
\begin{equation}
\label{EqApproxPexp}
\psi_{\gamma_{\theta_0,f_0}}(\rho)\asymp 
\begin{cases}
1, &\text{ if } \b+2q-1> \a+1/p,\\
(\log(1/\rho))^{(2-p)/2},& \text{ if } \b+2q-1=\a+1/p,\\
(1/\rho)^{\frac{p(\a +1/p-(\b+2q-1))}{\b+2q-1}},&\text{ if } \b+2q-1<\a+1/p.
\end{cases}
\end{equation}
The first case in this display is similar to the case of exponentially decreasing $g_k$. We apply Proposition~\ref{prop: GP} (i) 
and obtain the BvM theorem without needing a condition on the contraction rate. 
In the other cases we apply (ii) or (iii) of Proposition~\ref{prop: GP} and need to verify the condition $n\e_n^2\rho_n^2\to 0$ next
to the condition that $\psi_{\g_{\q_0,f_0}}(\rho_n)\ll n^{1/2}$ when $p=1$ or $\psi_{\g_{\q_0,f_0}}(\rho_n)\le n^{p/2}\rho_n^p$  when $p\in (1,2]$.
In the middle case of \eqref{EqApproxPexp} the latter conditions are verified for $\rho_n\asymp n^{-1/2}$ and $\rho_n\asymp n^{-1/2}(\log n)^{(2-p)/(2p)}$,
for $p=1$ and $p\in (1,2]$, respectively, and hence the condition $n\e_n^2\rho_n^2\to 0$ requires just $\e_n\to0$ or $\e_n\to0$ faster
than a logaritmic rate. More interesting is the third case in the display, which we consider separately when $p=1$ or $p\in (1,2]$.
In both cases $\a+1/p>\b+2q-1$ implies that $\a\ge\b$, since $2q-1-1/p\ge 1\ge 0$ by assumption,
so that $\e_n\asymp n^{-\b/(1+2\b+p(\a-\b))}$, by \eqref{EqRatepExp}.

In the case $p=1$ we apply (ii) of Proposition \ref{prop: pexp}, choosing  $\rho_n$ so that $\psi_{\gamma_{\theta_0,f_0}}(\rho_n)\ll n^{1/2}$. 
By \eqref{EqApproxPexp} this is satisfied in the case that $\a>\b+2q-2$ for 
$\rho_n\gg n^{-(\beta+2q - 1)/(2\a-2\b-4q+4)}$. By a long calculation the condition $n\eps_n^2\rho_n^2\to 0$ can then be seen
to be verified for $\alpha < 3\beta + 4q -3$. 

In the case $p\in(1,2]$ we apply (iii) of Proposition \ref{prop: pexp}, choosing  $\rho_n$ such that  $\psi_{\gamma_{\theta_0,f_0}}(\rho_n)\leq n^{p/2}\rho_n^p$.
By \eqref{EqApproxPexp} this is satisfied for $\rho_n\asymp n^{-(\beta+2q-1)/(2(\alpha+1/p))}$. 
The condition $n\eps_n^2\rho_n^2\to 0$ then translates into an upper bound on $\alpha$ involving $\beta, q$ and $p$. (It can be written
as a quadratic inequality in $\a$; we omit the details.)  If we write the upper bound as 
$\alpha < A(\beta,q,p)$, then one can show that for a fixed value of $p$, i.e.\ a fixed prior, the number $A(\beta,q,p)$ increases in $\beta$ and $q$. 
This behavior is to be expected for $\beta$, while for $q$ it suggests that more ill-posedness leads to a wider validity of the BvM phenomenon. 
For $p=2$ it can be calculated that $A(\beta,q,2)=2\beta + 2q -3/2$, which is identical to the bound obtained 
for the Gaussian series prior in Example \ref{ex: seriesprior}. One can show that, as $p\downarrow 1$, 
\[ A(\beta,q,p)\to \beta+q-\frac{3}{2} + \sqrt{2\beta^2+2\beta(2q-1)+(q-1/2)^2}.\]
Whenever $\beta, q$ are not small, this limit dominates $A(\beta,q,2)$. 
For example, if we only use that $q\geq 3$, then this happens if $\beta > \sqrt{27}/2-2$. 

The preceding calculations show that a prior with smaller $p$ will typically lead to a wider range of smoothness levels $\b$ of the true function $f_0$
and the prior parameter $\a$ for which the BvM phenomenon is valid.
\end{ex}
 
\section{Posterior simulations}
\label{sec:simulations}

In this section, we present some simulations to illustrate the application of our BvM results. To be brief, we will content ourselves with simulations for the heat equation example in Section~\ref{sec:applications} when putting a gaussian prior on $f$. That is, we will be demonstrating the result of Proposition~\ref{prop: GP}. 

We begin by reformulating the model to use the series representation for the prior, which is more convenient for sampling using Monte-Carlo Markov Chain (MCMC) methods. We select a ground truth $(\theta_0,f_0)$ and implement a Metropolis-Hastings (MH) algorithm to sample from the marginal posterior distribution of $\theta$. We repeat this process for different prior regularities $\alpha$ and examine the resulting marginal posterior distribution of $\theta$. We also provide trace plots to assess the convergence of our algorithm.

For a given pair of parameters, the signal-in-white noise model we consider is equivalent to observing the following samples for all $k\in\mathbb{N}$:

\begin{align*}
	X^{(1)}_k &= f_k + \frac{1}{\sqrt{n}}\zeta^1_k, \\
	X^{(2)}_k &= e^{-ck^2\theta}f_k + \frac{1}{\sqrt{n}}\zeta^2_k,
\end{align*}
where the $\zeta^1_k$ and $\zeta^2_k$ are i.i.d.\ standard Gaussians. In practice, we compute an approximation of the posterior by considering only the first $m$  observations. 
This yields the following expression for the approximate likelihood:
\begin{align*}
\ell(\theta,f\mid X^n)= \dfrac{n^{K}}{\sqrt{2\pi}^{2K}}\prod_{k\leq m} \exp\left(-\dfrac{n}{2}(X^{(1)}_k-f_k)^2 - \dfrac{n}{2}(X^{(2)}_k- e^{-ck^2\theta}f_k)^2\right).
\end{align*}
We put a uniform prior on $\theta$ and consider the prior series representation of the GP prior on $f$ already mentioned in Example \ref{ex: seriesprior}.
\begin{align*}
	\pi^\sigma_f\sim \sum_{k=1}^{\infty} \sigma_k\nu_ke_k,
\end{align*}
where the $\nu_k$ are i.i.d $N(0,1)$ random variables and where $\sigma_k = (k+1)^{-1/2-\alpha}$ where $\alpha>0$ is the chosen regularity of the prior. It follows that the prior is distributed as follows
\begin{align*}
	(\theta,f_1,\cdots,f_m) \sim U(\Theta)\times \prod_{k=1}^m N\left(0,(k+1)^{-1-2\alpha}\right).
\end{align*}
Furthermore, by integrating $\Pi(\theta,f)\ell(\theta,f|X^n)$ with respect to to the $f_k$'s, we obtain that the marginal posterior for $\theta, \Pi(\theta\mid X^n)$ is proportional to the following quantity:
\begin{align}\label{eq:post_theta}
	\prod_{k\leq m}(n+ne^{-ck^2\theta}+(k+1)^{2\alpha+1})^{-1/2} \exp\left(\dfrac{(nX^{(1)}_k+nX^{(2)}_ke^{-ck^2\theta})^2}{2(n+ne^{-ck^2\theta}+(k+1)^{2\alpha+1})}\right).
\end{align}
Our MH algorithm uses \eqref{eq:post_theta} as a target distribution. For the proposal step, we sample from a normal distribution centered at the previous proposal value with standard deviation $\sigma_\theta$. The details of our sampling procedure are described in Algorithm~\ref{alg:MH}.
\newline

\begin{algorithm}[]
	
	\SetAlgoLined
	\KwData{Observed data $\mathbf{X^{(1)},X^{(2)}}$, signal-to-noise ratio $\mathbf{n}$, regularity $\boldsymbol{\alpha}$, initial value $\boldsymbol{\theta^{(1)}}$, number of iterations $\mathbf{T}$, std.deviation for proposal distribution $\boldsymbol{\sigma_\theta}$. }
	\KwResult{Samples from the posterior distribution.}
	
	\For{$t=2$ to $\mathbf{T}$}{
		
		Sample $\theta^*$ from $N(\theta^{(t-1)},\boldsymbol{\sigma_\theta}^2)$
		
		Compute the acceptance ratio $A(\theta^*,\theta^{(t-1)})=\frac{\Pi(\theta^* | \mathbf{X^{(1)}},\mathbf{X^{(2)}})}{\Pi(\theta^{(t-1)} | \mathbf{X^{(1)}},\mathbf{X^{(2)}})}$ using \eqref{eq:post_theta}
		
		Sample $u$ from the uniform distribution over $(0,1)$
		
		\If{$u<A(\theta^*,\theta^{(t-1)})$}{
			$\theta^{(t)} = \theta^*$
		}
		\Else{
			$\theta^{(t)} = \theta^{(t-1)}$
		}
	}
	\caption{Metropolis-Hastings Algorithm}
	\label{alg:MH}
\end{algorithm}

We fix $T=1$, let $\theta_0=0.01$ and consider the function $f_0$ defined as follows:

\begin{align*}
	f_0(x)=\sum_{k=1}^\infty k^{-2}e_k(x).
\end{align*}

\begin{figure}
	\centering
	\includegraphics[width=\textwidth, height=4cm]{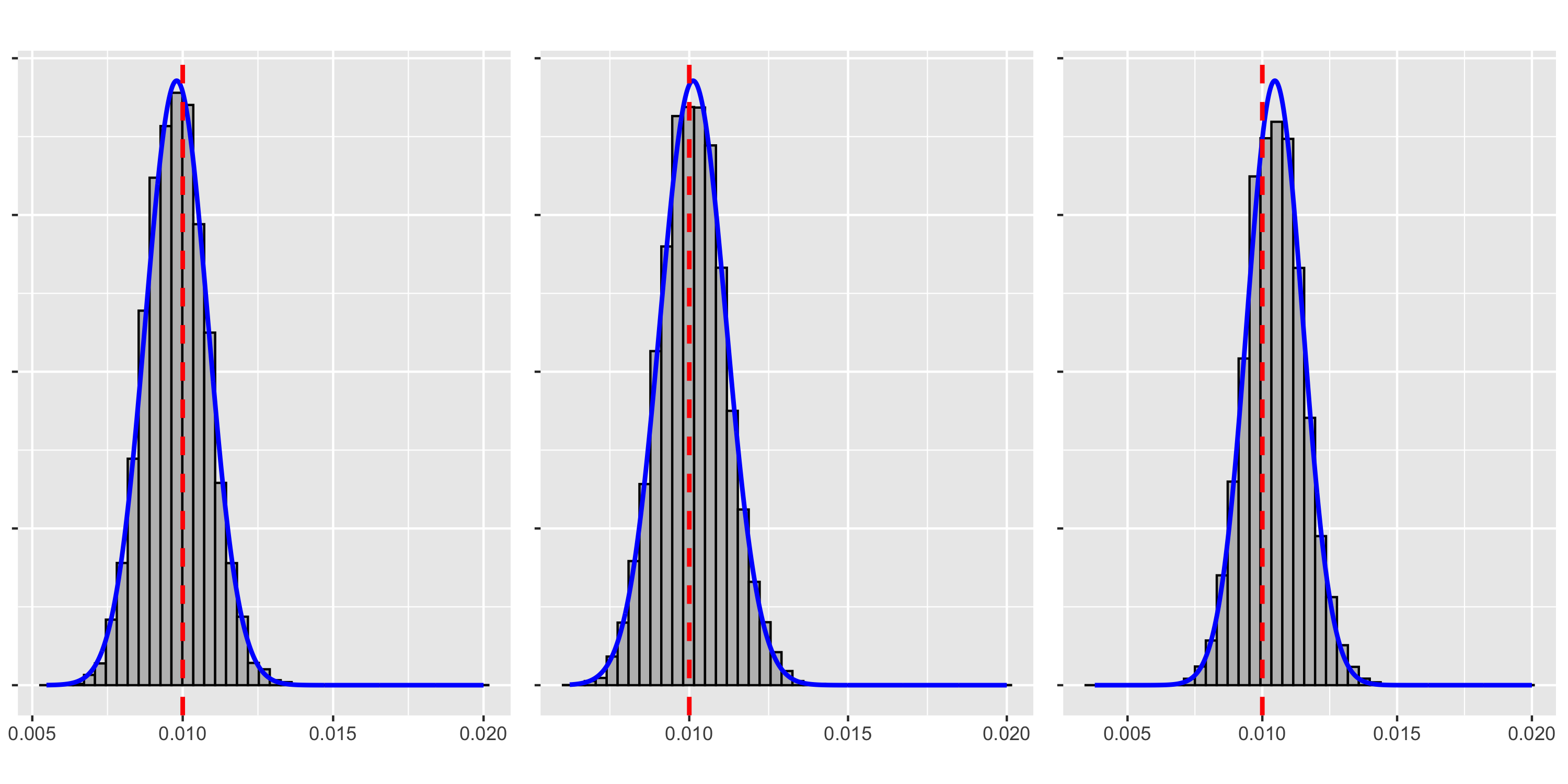}
	\caption{Approximations of the marginal posterior of $\theta$ obtained when running algorithm \ref{alg:MH} on three different datasets with $\alpha=1$. The red line marks the true value of $\theta$ while the blue curve is the theoretical limiting distribution in Theorem~\ref{thm: main}. That is, a normal distribution centered at an efficient estimator (here the posterior mean) with variance the efficient Fisher information.}
	\label{fig:goodAlpha}
\end{figure}

\begin{figure}
	\centering
	\includegraphics[width=\textwidth,height=4cm]{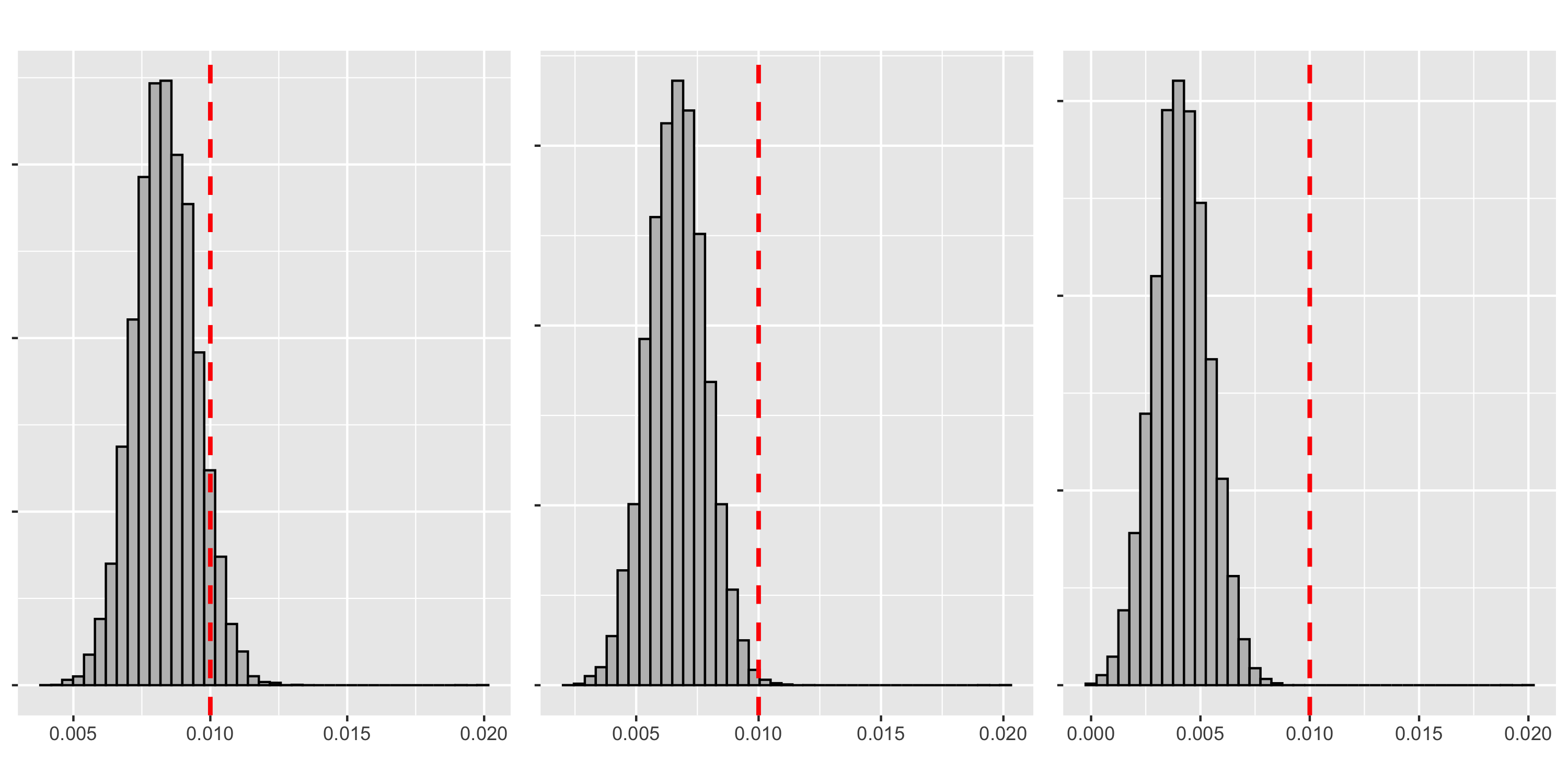}
	\caption{Approximations of the marginal posterior of $\theta$ for priors with regularity $\alpha$ equal to (from left to right) 2.6, 3.0 and 3.4. All three values of $\alpha$ are outside the BvM zone predicted by Theorem~\ref{thm: main} The approximations are realized for the same dataset.}
	\label{fig:ALPHAS_outisde}
\end{figure}

We have that $f$ is Sobolev smooth of order $\beta = 3/2$. Our BvM theorem is guaranteed for $0.5< \alpha < 2.5$. We fix $m=100$ and  generate data vectors $X^{(1)}$ and $X^{(2)}$ with signal-to-noise ratio $n=10^5$.  We run Algorithm \ref{alg:MH} for $10^5$ iterations with a burn in period of 1000. For a correctly chosen $\alpha$ in the zone ($\alpha=1$), we run Algorithm \ref{alg:MH} on different datasets and observe that the marginal posterior of $\theta$ indeed seems to satisfy the BvM (Figure \ref{fig:goodAlpha}).  We also plot the resulting marginal posterior distributions for $\theta$ for different values of $\alpha$ outside of the zone prescribed by the BvM theorem in Figure \ref{fig:ALPHAS_outisde}, that is for $\alpha>2\beta-1/2$. As anticipated, we observe the appearance of a bias which becomes larger when $\alpha$ increases. This process was repeated with different datasets to confirm the observations.

Overall, the presented figures validate the predictions of our semi-parametric BvM. They also exhibit the positive relationship between the prior regularity $\alpha$ and the magnitude of a bias in the marginal posterior distributions which can be observed for $\alpha$ outside the BvM zone. 


To verify the convergence of our sampling algorithm, we present two of the trace plots of the chains of $\theta$ obtained during the sampling process for different values of $\alpha$ and on different datasets. These trace plots showcase no particular behavior, and no burn in period seems to really be required to obtain convergence. This contributes to the evidence that our algorithm is indeed converging and sampling from the posterior distribution. 

\begin{figure}
	\centering
	\includegraphics[width=0.8\textwidth,height=4cm]{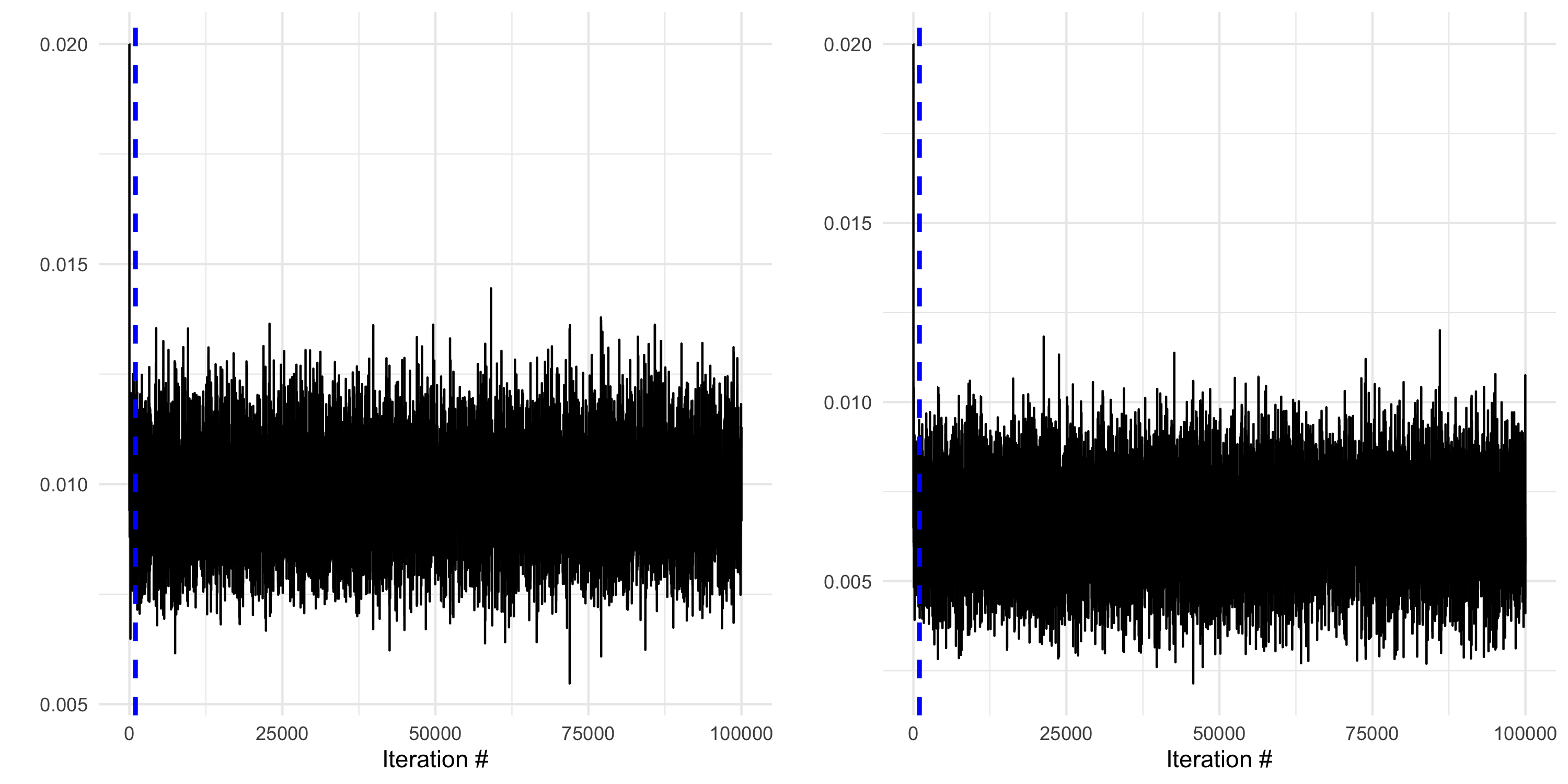}
	\caption{Trace plots of $\theta$ for two MH algorithm runs. The blue lines mark the burn in period selected.}
	\label{fig:traceplots}
\end{figure}

\section{Proofs}\label{sec: proof}

\subsection{Proof of Lemma~\ref{lem:LAN}}
Two applications of formula \eqref{EqLikelihood} for the log-likelihood give, for $(\tau,h)\in\Theta\times H$, 
\begin{align*}
	\begin{split}
		\log{\frac{dP_{\tau,h}^n}{dP_{\theta,f}^n}} &= \sqrt{n}\ip{h-f, \dot W^{(1)}} - \frac{n}{2}\|h-f\|^2 \\
		&\qquad +\sqrt{n}\ip{K_\tau h - K_\theta f, \dot W^{(2)}} - \frac{n}{2}\|K_\tau h-K_\theta f\|^2.
	\end{split}
\end{align*}
Therefore, the left side of the lemma can be written as 
\begin{align*}
&t\ip{g , \dot W^{(1)}} + \langle K_{\theta+t/\sqrt{n}}(f+tg/\sqrt{n})-K_\theta f, \dot W^{(2)}\rangle \\
&\qquad\qquad - \frac{t^2}{2}\|g\|^2- \frac{1}{2}\|K_{\theta+t/\sqrt{n}}(f+tg/\sqrt{n})-K_\theta f\|^2
\end{align*}
By assumption $s^{-1}\bigl(K_{\theta+as} (f+sg) - K_{\theta} f\bigr) \to a\dot K_\theta f + K_\theta g$ in $H$, as $s\to 0$.
This allows to expand the second and fourth terms in the preceding display further, as in the lemma, where
we use that $\ip{h,\dot W^{(2)}}$ is normally distributed with mean zero and variance $\|h\|^2$.

The Fisher information $M(g):=\|g\|^2 + \|\dot K_\theta f + K_\theta g\|^2$ satisfies, for any $g,h\in H$,
\[
M(g+h)-M(g) = \|h\|^2+\|K_\theta h\|^2 +2\ip{g,h}-2\ip{\dot K_\theta f - K_\theta g, K_\theta h}.
\]
This is clearly nonnegative if $\ip{g,h}-\ip{\dot K_\theta f - K_\theta g, K_\theta h} = 0$. For $g=-\gamma_{\theta,f}$ the latter
is true for every $h\in H$, whence $g$ minimizes $M$ over $H$.

The minimal value is clearly positive if $\g_{\q,f}\not=0$. If $\g_{\q,f}=0$, then $K_\q\g_{\q,f}=0$ 
and the minimal value is  $\|\dot K_\theta f \|^2>0$, which is positive if $\dot K_\q f\not=0$.

\subsection{Proof of Theorem~\ref{thm: main}}

The theorem follows from Theorem~12.9 in \cite{vaartghosal}, which is an adaptation of results by \cite{castillo}.
We apply Theorem~12.9 with the least favorable transformation $(\theta,f)\mapsto (\theta_0, f+(\theta-\theta_0)\gamma_n)$.
Then its condition (12.14) is satisfied in view of our assumption \eqref{eq: priorshift}, the posterior consistency
conditions in Theorem~12.9 are also copied in our conditions, and we need only verify (12.13) in \cite{vaartghosal}. 
Here we may assume that the sets $\Theta_n\times H_n$ are contained in shrinking balls 
$\{(\q,h)\in \Theta\times H: \|\q-\q_0\|<\e_n,\|h-f_0\|<\e_n\}$ of $(\q_0,f_0)$, for some $\e_n\to0$,
by our assumption of posterior consistency.

Straightforward computations using \eqref{EqLikelihood} yield
\[
\log \frac{{dP_{\theta,f}^n}}{dP_{\theta_0, f+(\theta-\theta_0)\gamma_n}^n} 
= \sqrt{n}(\theta-\theta_0) G_{\theta_0}(f,\gamma_n) - \frac{n}{2}|\theta-\theta_0|^2 \Tilde I_{\theta_0, f}(\gamma_n) + R_n(\theta,f),
\]
where 
\begin{align*}
G_{\theta_0}(f,g)&= \ip{\dot W^{(1)}, -g} + \ip{\dot W^{(2)}, \dot K_{\theta_0}f - K_{\theta_0} g},\\
\Tilde I_{\theta_0, f}(g)&= \|\dot K_{\theta_0}f\|^2 - \|K_{\theta_0}g\|^2 - \|g\|^2,\\
R_n(\theta,f) &= \sqrt{n}\ip{\dot W^{(2)}, (K_{\theta}-K_{\theta_0} - (\theta-\theta_0)\dot K_{\theta_0})f}\\
		&\qquad + n(\theta-\theta_0)\ip{f-f_0, (I+K_{\theta_0}^*K_{\theta_0})(\gamma_n-\gamma_{\theta_0,f})} \\
		&\qquad- n\ip{(K_{\theta}-K_{\theta_0} - (\theta-\theta_0)\dot K_{\theta_0})f, K_{\theta_0}(f-f_0)}\\
		&\qquad +\frac{n}{2}(\theta-\theta_0)^2\|\dot K_{\theta_0}f\|^2- \frac{n}{2}\|(K_{\theta}-K_{\theta_0})f\|^2.
\end{align*}
Therefore, condition (12.14) of Theorem~12.9 in \cite{vaartghosal} is satisfied, with $\tilde G_n:=G_{\theta_0}(f_0,\gamma_{\theta_0,f_0})$
and $\tilde I_n=\tilde I_{\q_0,f_0}$, if 
\begin{align}
\sup_{(\theta,f)\in \Theta_n\times H_n} \frac{\sqrt{n}(\theta-\theta_0)| G_{\theta_0}(f, \gamma_n)-G_{\theta_0}(f_0,\gamma_{\theta_0,f_0})|}{1+n(\theta-\theta_0)^2} &\gaatp 0,\label{eq: LANapprox}\\
\sup_{(\theta,f)\in \Theta_n\times H_n} \frac{n(\theta-\theta_0)^2|\Tilde I_{\theta_0,f} (\gamma_n)-\Tilde I_{\theta_0, f_0}|}{1+n(\theta-\theta_0)^2} &\to 0,\label{EqContFI}\\
\sup_{(\theta,f)\in \Theta_n\times H_n} \frac{R_n(\theta,f)}{1+n(\theta-\theta_0)^2} &\gaatp 0.\label{eq: remaindercondition}
\end{align}
The Bernstein-von Mises theorem is then satisfied at $(\q_0,f_0)$ with 
$\Delta_{\theta_0, f_0}^n:= \Tilde I_{\theta_0,f_0}^{-1}G_{\theta_0}(f_0, \gamma_{\theta_0,f_0})$, which is exactly $N(0, \Tilde I_{\theta_0,f_0}^{-1})$
distributed. (Note that the statement of Theorem~12.9 in \cite{vaartghosal} needs an extra factor $n^{-1/2}$, 
as acknowledged in the list of errata \cite{vaartghosal_errata}.)

In \eqref{eq: LANapprox} we bound the quotient $\sqrt n|\q-\q_0|/(1+n|\q-\q_0|^2)$ by 1, and split 
the variable $ G_{\theta_0}(f, \gamma_n)-G_{\theta_0}(f_0,\gamma_{\theta_0,f_0}) $ 
into the two terms 
$T_{1,n}:=\ip{\dot W^{(1)}, \gamma_{\theta_0,f_0}-\gamma_n} + \ip{\dot W^{(2)}, K_{\theta_0}(\gamma_{\theta_0,f_0}-\gamma_n)}$
and $T_{2,n}(f):=\ip{\dot W^{(2)}, \dot K_{\theta_0}(f-f_0)}$. The first term $T_{1,n}$ is a Gaussian variable with mean zero and
variance tending to zero by the assumption that $\g_n\to\g_{\q_0,f_0}$ and Assumption~\ref{asspt:main}(i),
whence it converges to zero in probability. The second term $T_{2,n}(f)$ is a centered Gaussian process
 indexed by the set $\mathcal{T}_n= \{\dot K_{\theta_0}(f-f_0): f\in H_n\}$ with square intrinsic metric
$\var\bigl(T_{n,2}(f)-T_{n,2}(g)\bigr)=\|\dot K_{\theta_0}(f-g)\|^2$.
 bounded above by the square $H$-metric, by Assumption~\ref{asspt:main} (i). 
By  Assumption~\ref{asspt:main} (ii), (iii), the $S$-norm of the function  $\dot K_{\theta_0}(f-f_0)$ is bounded by a multiple of $\|f-f_0\|$,
and the same is true for the $H$-norm. Thus $\mathcal{T}_n$ is contained in a multiple of the unit ball of $S$
and its $H$-diameter is bounded by $\e_n$ for some $\e_n\ra 0$, by the posterior consistency assumption.
It follows, by Dudley's bound (\cite{Dudley73}, or Corollary~2.2.9 in \cite{vdV_wellner2023}), that
\[
\mathbb{E}\sup_{f\in \mathcal{T}_n}|T_{2,n}(f)|\lesssim \int_0^{\e_n} \sqrt{\log N(\eps, \{h: \|h\|_S\le 1\}, \|\cdot\|)}\,d\eps.
\] 
The integral tends to zero by assumption \eqref{eq: integralCondition}.
Thus \eqref{eq: LANapprox} is verified.

For \eqref{EqContFI} it suffices to show that $\tilde I_{\q_0,f_0}-\Tilde I_{\theta_0, f}(\gamma_n)$ tends to
zero, uniformly in $f\in H_n$. Because $\tilde I_{\q_0,f_0}=\Tilde I_{\theta_0, f_0}(\gamma_{\theta_0,f_0})$, and $\dot K_{\q_0}$
and $K_{\q_0}$ are continuous by Assumption~\ref{asspt:main}, this can be reduced to 
$\sup_{f\in H_n}\|f-f_0\|\to 0$ and $\|\gamma_n-\gamma_{\theta_0,f_0}\|\to 0$, which are true by
the posterior consistency assumption.

Finally we prove \eqref{eq: remaindercondition}, where we consider the four terms defining $R_n(\q,f)$ separately.

For the first term we bound the quotient $\sqrt n/(1+n|\q-\q_0|^2)$ by  $1/|\q-\q_0|$, and next consider  the class of functions
\[
\GG_n = \Bigl\{ g_{\q,f}:=\frac{(K_\theta-K_{\theta_0}-(\theta-\theta_0)\dot K_{\theta_0}) f}{\theta-\theta_0}: 
(\theta, f) \in \Theta_n\times H_n \Bigr\}.
\]
By Assumption~\ref{asspt:main} (iii) the $S$-norm of the function $g_{\q,f}\in\GG_n$ is bounded above by a 
multiple of $|\q-\q_0|\,\|f\|$, and hence so is the $H$-norm. By the assumption of posterior contraction
this is bounded above by a multiple of $\e_n$ uniformly in $g_{\q,f}\in\GG_n$, for some $\e_n\ra 0$.  The intrinsic
distance of the process $\langle \dot W^{(2)},g\rangle$ is given by the $H$-norm of $\|g\|$.
By a similar argument as previously, using Dudley's bound, it therefore follows that $\sup_{g\in\GG_n}|\ip{\dot W^{(2)},g}|$
tends to zero in probability. 

For the second term, we bound $n|\q-\q_0|/(1+n|\q-\q_0|^2)$ by $\sqrt n$,
and thus reduce it to assumption \eqref{eq: lfdrate}.

By the  Cauchy-Schwartz inequality and item (iii) of Assumption~\ref{asspt:main}, the third term 
is bounded above by a multiple of $n(\theta-\theta_0)^2\|f\|\,\|f-f_0\|$, which is $o(1+n(\theta-\theta_0)^2)$,  uniformly in $f\in H_n$. 

Finally, by item (ii) of Assumption~\ref{asspt:main}, the fourth term is bounded above by a multiple of $n|\theta-\theta_0|^3\|f\|$, 
which is $o(1+n(\theta-\theta_0)^2)$ uniformly on $\Theta_n\times H_n$. 

\subsection{Proof of Lemma~\ref{thm: main2}}

The ``intrinsic''  (square) metric for the  (extended) white noise model \eqref{eq: obs1}--\eqref{eq: obs2} is
$$d^2\bigl((\theta_1, f_1), (\theta_2, f_2)\bigr) = \|f_1-f_2\|^2 + \|K_{\theta_1} f_1 - K_{\theta_2}f_2\|^2.$$
We show below that under (ii) and (iii) of Assumption~\ref{asspt:main},
there exists a constant $A>0$ such that, for every $\q\in\Theta$ and $f\in H$,
\begin{equation}
|\theta-\theta_0| + \|f-f_0\| \le A\, d\bigl((\theta, f), (\theta_0, f_0)\bigr).
\label{EqComparisonMetrics}
\end{equation}
It follows that a rate of contraction for $d$ implies the same rate of contraction for the natural metric,
given on the left. Such a rate of contraction $\e_n$ follows by a straightforward adaptation of 
Theorem~8.31 of \cite{vaartghosal}, under the conditions $n\e_n^2\ge 1$ and the existence
of sets $\FF_n\subset H$ such that 
\begin{align}
	\label{eq: cc1}
\Pi\bigl((\theta, f): d((\theta, f), (\theta_0, f_0)) \le \eps_n\bigr) & \ge e^{- n\eps^2_n/64},\\
\label{eq: cc2} 
\sup_{\e>\eps_n}\log N\bigl(\eps/8, \{(\q,f)\in\Theta \times \FF_n: d((\theta, f), (\theta_0, f_0)) \le \eps\}, d\bigr) & \le n\eps^2_n,\\
	\label{eq: cc3} \Pi\bigl((\theta, f) \not \in \Theta\times\FF_n\bigr) & \le e^{-n\eps^2_n}.
\end{align}
The last condition \eqref{eq: cc3} is trivially implied by \eqref{eq: thmPrior2}. 

For the verification of \eqref{eq: cc1} and \eqref{eq: cc2}, we use that the intrinsic distance 
$d$ satisfies, for all $\theta_1, \theta_2 \in \Theta$ and $f_1, f_2 \in \FF_n$,
by (iv) of assumption Assumption~\ref{asspt:main},
\[
d\bigl((\theta_1, f_1), (\theta_2, f_2)\bigr) \le D_3  |\theta - \theta_0|(\|f_1\|\wedge \|f_2\|) +(D_3+1) \|f_1-f_2\|.
\]
Thus the set in \eqref{eq: cc1} contains the set of all $(\q,f)$ with $|\q-\q_0|\,\|f_0\|<c\e_n$ and $\|f-f_0\|<c\e_n$,
for some $c=1/(2D_3+2)$. Since the prior density for $\theta$ is assumed to be bounded away from $0$, 
\eqref{eq: cc1} is implied by \eqref{eq: thmPrior}. 

Since by assumption the elements of $\FF_n$ possess norm bounded by $R_n= e^{n\e_n^2/4}$,
on $\Theta\times\FF_n$ the metric $d$ is bounded above by a 
multiple of $ |\theta - \theta_0|R_n +\|f_1-f_2\|$. It follows that, for some constant $c$,
\[
\log N(\eps_n/c, \Theta \times \FF_n, d) \le \log N(\eps_n, \Theta, R_n|\cdot|)
+ \log N(\eps_n, \FF_n, \|\cdot\|). 
\]
The first term on the right is of the order $\log (R_n/\eps_n)\le n\eps_n^2$ by the assumption on $R_n$.
The second term on the right is bounded by a constant times 
$n\eps^2_n$ by \eqref{eq: thmCover}. 

We finish the proof of the first assertion of the lemma by deriving \eqref{EqComparisonMetrics}.
We can write $K_\q f-K_{\q_0}f_0=(\q-\q_0)\dot K_{\q_0}f_0+K_{\q_0}(f-f_0)+R(\q,f)$,
for the remainder $R(\q,f)=(K_\q-K_{\q_0})f_0-(\q-\q_0)\dot K_{\q_0}f_0+(K_\q-K_{\q_0})(f-f_0)$.
By Assumption~\ref{asspt:main} (ii), (iii), we have $\|R(\q,f)\|\le c (\q-\q_0)^2+c|\q-\q_0|\,\|f-f_0\|$, where the
proportionality constant $c$ depends on $f_0$. We conclude that, for $\q\not=\q_0$,
\begin{align*}&\|K_\q f-K_{\q_0}f_0\|+\|f-f_0\|\\
&\qquad\ge |\q-\q_0|\Bigl(\Bigl\|\dot K_{\q_0}f_0+K_{\q_0}\frac{f-f_0}{\q-\q_0}\Bigr\|
+\Bigl\|\frac{f-f_0}{\q-\q_0}\Bigr\|-c|\q-\q_0|-c\|f-f_0\|\Bigr).
\end{align*}
Here for any $g$,
$$\|\dot K_{\q_0}f_0+K_{\q_0}g\|+\|g\|\ge \sqrt{\|\dot K_{\q_0}f_0+K_{\q_0}g\|^2+\|g\|^2}
\ge \sqrt{\tilde I_{\q_0,f_0}},$$
which is strictly positive, by Lemma~\ref{lem:LAN}.
Thus the right side of the second last display is bounded below
by $|\q-\q_0|\bigl(\sqrt{\tilde I_{\q_0,f_0}}-c \eta\bigr)$, if $|\q-\q_0|+\|f-f_0\|<\eta$.
This proves \eqref{EqComparisonMetrics} for $(\q, f)$ in a neighborhood of $(\q_0,f_0)$.
Failure of the inequality outside the neighborhood would entail existence of a sequence
$(\q_m,f_m)$ with $\bigl(\|K_{\q_m} f_m-K_{\q_0}f_0\|+\|f_m-f_0\|\bigr)/\bigl(|\q_m-\q_0|+\|f_m-f_0\|\bigr)\to0$.
By compactness of $\Theta$ we can assume that $\q_m\to\q_1\not=\q_0$. It follows that $K_{\q_m}f_m\ra K_{\q_0}f_0$ and
$f_m\ra f_0$ and hence $K_{\q_1}f_0=K_{\q_0}f_0$, which is excluded by assumption.

The second assertion of the lemma, on the posterior distribution for known value of $\q_0$, can be proved
similarly, omitting the $\q$-part of the arguments.

\subsection{Proof of Lemma \ref{lemma: priorshift}}

By Proposition~2.7 in \cite{pexp}, for any $h\in\mathcal{Q}$ the measure $\pi_{f+h}$ is absolutely continuous relative to $\pi_f$ with
Radon-Nikodym derivative satisfying 
\begin{equation*}
\log\frac{d\pi_{f+h}}{d\pi_f}(f) = \frac{1}{p}\sum_{k=1}^\infty \Bigl(\left|\frac{f_k}{\sigma_k}\right|^p - \left|\frac{f_k-h_k}{\sigma_k}\right|^p\Bigr).
\end{equation*}
Two applications of this formula give, for $\g_n\in \mathcal{Q}$,
\begin{equation}
\label{EqRNpexp}
\log\frac{d\pi_{f+(s+t)\g_n}}{d\pi_{f+s\g_n}}(f) 
= \frac{1}{p}\sum_{k=1}^\infty \Bigl(\left|\frac{f_k-(s+t)\g_{n,k}}{\sigma_k}\right|^p - \left|\frac{f_k-s\g_{n,k}}{\sigma_k}\right|^p\Bigr).
\end{equation}

\subsubsection*{\textbf{Case $p=1$.}} 
For $p=1$ we can apply the triangle inequality to see that the right side of \eqref{EqRNpexp} is bounded in absolute value by
$t\sum_{k=1}^\infty|\g_{n,k}/\sigma_k|= t\|\g_n\|_{\mathcal{Z}}$.
It follows that 
\begin{align*}
\frac{\bigl|\log\bigl(d\pi_{f+(s+t)\g_n}/d\pi_{f+s\g_n})\bigr(f)\bigr|}{1+nt^2}
&\leq \frac{|t|\|\gamma_n\|_{\mathcal{Z}}}{1+nt^2} \leq \frac{\|\gamma_n\|_{\mathcal{Z}}.}{\sqrt n}.
\end{align*}
Condition \eqref{eq: priorshift} is thus satisfied (surely) if $n^{-1/2}\|\gamma_n\|_{\mathcal{Z}}\to 0$.

\subsubsection*{\textbf{Case $p\in(1,2]$.}}  
By \eqref{EqRNpexp}, Lemma~\ref{claim1} applied with  $x=f_k/\sigma_k$, $\eps=t\g_{n,k}/\sigma_k$ and $\d=s\g_{n,k}/\sigma_k$,
and the triangle inequality,
\begin{align*}
\left|\log \frac {d\pi_{f+(s+t)\g_n}}{d\pi_{f+s\g_n}}(f)\right|
&\le  \frac{1}{p}\Bigl|\sum_{k=1}^\infty \sign (f_k)\left|\frac{f_k}{\sigma_k}\right|^{p-1}\frac{t\g_{n,k}}{\sigma_k}\Bigr|\\
&\qquad+2\sum_{k=1}^\infty\Bigl|\frac{t\g_{n,k}}{\sigma_k}\Bigr|^p
+2\sum_{k=1}^\infty\Bigl|\frac{s\g_{n,k}}{\sigma_k}\Bigr|^{p-1}\Bigl|\frac{t\g_{n,k}}{\sigma_k}\Bigr|\\
&= \frac{|t|}{p}\bigl|\sum_{k=1}^\infty W_k\bigr|+2\bigl(|t|^p+|s|^{p-1}|t|\bigr)\|\g_n\|^p_{\mathcal{Z}},
\end{align*}
for $W_k:=\sign (f_k)|f_k/\sigma_k|^{p-1}\g_{n,k}/\sigma_k$. 
It follows that the left side divided by $1+nt^2$  is bounded above by a multiple of 
$$\frac{\bigl|\sum_{k=1}^\infty W_k\bigr|}{\sqrt n}+\frac{\|\g_n\|^p_{\mathcal{Z}}}{\sqrt n^p}
+\frac{|s|^{p-1}\|\g_n\|^p_{\mathcal{Z}}}{\sqrt n}
\le K\eps_n\|\gamma_n\|_{\mathcal{Q}}+\frac{\|\gamma_n\|_{\mathcal{Z}}^p }{n^{p/2}}
+\frac{\e_n^{p-1}\|\gamma_n\|_{\mathcal{Z}}^p }{\sqrt n},$$
on the event $V_n := \bigl\{f\in H: \bigl|\sum_{k=1}^\infty W_k \bigr|\leq K\sqrt{n}\eps_n\|\gamma_n\|_{\mathcal{Q}} \bigr\}$
and for $|s|< \e_n$. Every of the three terms on the right tends to zero if the relations in (ii) of the lemma are true.
Furthermore, since $\|\cdot\|_{\mathcal{Q}}\le \|\cdot\|_{\mathcal{Z}}$, the condition $\eps_n\|\gamma_n\|_{\mathcal{Z}}\to0$
implies that the first term tends to zero. It also implies that the other terms tend to zero, as $\e_n\ge 1/\sqrt n$ and
the third term can be factorized as $(\e_n\|\gamma_n\|_{\mathcal{Z}})^{p-1} \, (\|\gamma_n\|_{\mathcal{Z}}/\sqrt n)$.

To conclude the proof we show that the posterior mass of the $V_n$'s tends to 1. 
By Lemma~\ref{claim2} the variable $\sign(Z_k)|Z_k|^{p-1}$, where $f_k=\s_k Z_k$, is sub-Gaussian.
Because these variables are centered and independent, it follows that $\sum_{k=1}^\infty W_k=\sum_{k=1}^\infty\sign(Z_k)|Z_k|^{p-1}\g_{n,k}/\s_k$
 is sub-Gaussian as well, with variance proxy a multiple of $\sum_k (\g_{n,k}/\sigma_k)^2=\|\g_n\|_{\mathcal{Q}}^2$ 
(e.g. \cite{vdV_wellner2023}, Proposition A.1.6). Therefore 
$\pi_f(V_n^c)=\Pr\bigl(\bigl|\sum_{k=1}^\infty W_k \bigr|\ge K\sqrt{n}\eps_n\|\gamma_n\|_{\mathcal{Q}} \bigr)\le e^{-3n\e_n^2}$, 
for sufficiently large $K$. Combined with \eqref{eq: thmPrior}, the remaining mass principle (Theorem~8.20 in \cite{vaartghosal}) gives
that $\Pi(\Theta\times V_n|X^n)\to 1$ in $P_{\theta_0,f_0}$ probability.

\begin{lemma}\label{claim1}
If $p\in (1,2]$, then for any  $x,\eps,\d \in\mathbb{R}$,
\begin{align*}
	\bigl| |x|^p - |x-\eps|^p - p\sign(x)|x|^{p-1}\eps\bigr| &\le 2p|\eps|^p,\\
	 \bigl| |x-\d|^p - |x-\d-\eps|^p - p\sign(x)|x|^{p-1}\eps\bigr| &\le 2p|\eps|^p+2p|\d|^{p-1}|\e|.
\end{align*}
\end{lemma}

\begin{proof}
For $\e<0$ the first inequality can be obtained from the inequality for $-\e$ and $-x$. Therefore,
without loss of generality, assume $\eps>0$. If  $x>\eps\ge 0$, then the left side is equal to the absolute value of
\[x^p - (x-\eps)^p-px^{p-1}\e = \int_{x-\eps}^x (ps^{p-1}-px^{p-1})\,ds.\]
Since $p\in(1,2]$ the integrand can be bounded as $p|s^{p-1}-x^{p-1}|\leq p|s-x|^{p-1}$, and it follows that the absolute value of the
right side is bounded above by $\e^p$. For $x<0<\e$, the left side of lemma is the absolute value of, for $y=-x$,
\[y^p-(y+\eps)^p +py^{p-1}\e= -\int_y^{y+\eps}(ps^{p-1}-py^{p-1})\,ds.\]
This is bounded in absolute value by $\e^p$ as before.
Finally if $0\le x\le \eps$, then $\bigl||x|^p-|x-\e|^p\bigr|\le \e^p$ and $\big|p\sign(x)|x|^{p-1}\e\bigr|\le p\e^p$ 
and hence the inequality is true in view of the triangle inequality.

To prove the second inequality, we start by applying the first inequality with $x-\d$ instead of $x$. In view of the triangle inequality
it next suffices to prove that $\bigl|p\sign(x-\d)|x-\d|^{p-1}\eps-p\sign(x)|x|^{p-1}\eps\bigr|\le 2p|\d|^{p-1}|\e|$, or
 $\bigl|\sign(x-\d)|x-\d|^{p-1}-\sign(x)|x|^{p-1}\bigr|\le 2|\d|^{p-1}$. For $\d<0$, this inequality can be obtained from
the inequality with $-\d$ and $-x$. Therefore, again assume without loss of generality that $\d>0$. If $x>\d>0$, then
the inequality is $|(x-\d)^{p-1}-x^{p-1}|\le 2\d^{p-1}$, which is true (without 2) because $0\le p-1\le1$. For $0\le x<\d$, the inequality
is $|-(\d-x)^{p-1}-x^{p-1}|\le 2\d^{p-1}$, which is true by the triangle inequality. For $x<0<\d$, the inequality is
$|-(\d+|x|)^{p-1}+|x|^{p-1}|\le 2\d^{p-1}$, which is true because $0\le p-1\le1$.
\end{proof}

\begin{lemma}
\label{claim2}
Let $1< p\le 2$. If $Z$ has density proportional to $z\mapsto e^{-|z|^p/p}$, then
the variable $X=\sign(Z)|Z|^{p-1}$ is sub-Gaussian with variance proxy $p$, i.e.\ for all $t>1$,
$$\mathbb{P}(|X|>t) \lesssim \frac{2e^{-t^2/p}}{t}.$$
\end{lemma}

\begin{proof}
The probability on the left is equal to $2\mathbb{P}(Z>s)$, for $s:=t^{\frac{1}{p-1}}$.
For  $s>0, p>1$, we have $x^{p-1}>s^{p-1}$ for all $x>s$ and hence
\[
\mathbb{P}(Z>s) \propto \int_s^\infty e^{-x^p/p}\,dx \le \frac{1}{s^{p-1}}\int_s^\infty x^{p-1}e^{-x^p/p}\,dx 
= \left[-e^{-x^p/p}\right]_s^\infty=\frac{e^{-s^p/p}}{s^{p-1}}.
\]
Plugging in the value of $s$ gives the bound $2e^{-t^{p/(p-1)}/p}/t$ on the probability in the lemma.
As $p\in(1,2]$, we have $p > 2(p-1)$, and we obtain the bound of the lemma.
\end{proof}

\begin{acks}[Acknowledgements]
	
	The authors would like to thank Frank van der Meulen for his helpful assistance and comments, as well as the two anonymous referees and the Associate Editor for their constructive comments that improved the quality of this paper.
	
\end{acks}

\begin{funding}
	This research is partly supported by NWO grant 613.009.134 awarded to J.H. van Zanten
	and the NWO Spinoza prize awarded to A.W. van der Vaart by the Netherlands Organisation for Scientific Research (NWO).
\end{funding}

\bibliographystyle{imsart-number.bst}
\bibliography{references.bib}

\end{document}